\documentclass[arma,notitlepage]{svjour}
\usepackage{
amsmath,upref,
amssymb,
amsxtra,amscd,enumerate
}
\usepackage[cp1250]{inputenc}
\usepackage{enumerate}
\usepackage{color}
\usepackage[mathscr]{eucal}
\newcounter{cnstcnt}
\newcommand{\novak}{%
\refstepcounter{cnstcnt}%
\ensuremath{c_{\thecnstcnt}}}
\newcommand{\starak}[1]{\ensuremath{c_{\ref{#1}}}}
\newcommand{\bfxi}{\mbox{\boldmath $\xi$}}
\newcommand{\bfeta}{\mbox{\boldmath $\eta$}}
\newcommand{\bfom}{\mbox{\boldmath $\varpi$}}

\newtheorem{lemm}{Lemma}
\newtheorem{Def}{Definition}
\newtheorem{prop}{Proposition}
\newtheorem{theo}{Theorem}

\newtheorem{rem}{Remark}

\newcommand{\QED}{\par\hfill$\square$\par}

\newcommand{\Bd}{\begin{Def}}

\newcommand{\Ed}{\end{Def}}

\newcommand{\Bp}{\begin{prop}}

\newcommand{\Ep}{\end{prop}}

\newcommand{\Br}{\begin{rem}}

\newcommand{\Er}{\end{rem}}

\newcommand{\Bt}{\begin{theo}}

\newcommand{\Et}{\end{theo}}

\newcommand{\Bl}{\begin{lemm}}

\newcommand{\El}{\end{lemm}}

\newcommand{\diver}{{\rm div}\,}

\newcommand{\bfx}{\mbox{\boldmath $x$}}

\newcommand{\bfell}{\mbox{\boldmath $\ell$}}

\newcommand{\bfy}{\mbox{\boldmath $y$}}

\newcommand{\bfphi}{\mbox{\boldmath $\varphi$}}

\newcommand{\bfpsi}{\mbox{\boldmath $\psi$}}

\newcommand{\bfv}{{\mbox{\boldmath $v$}} }

\newcommand{\bfn}{{\mbox{\boldmath $n$}} }

\newcommand{\bfu}{{\mbox{\boldmath $u$}} }

\newcommand{\bfq}{{\mbox{\boldmath $q$}} }

\newcommand{\bfw}{{\mbox{\boldmath $w$}} }

\newcommand{\bfa}{{\mbox{\boldmath $a$}} }

\newcommand{\bfg}{{\mbox{\boldmath $g$}} }

\newcommand{\bfr}{{\mbox{\boldmath $r$}} }

\newcommand{\bfk}{{\mbox{\boldmath $k$}} }

\newcommand{\bfV}{{\mbox{\boldmath $V$}} }

\newcommand{\bfW}{{\mbox{\boldmath $W$}} }

\newcommand{\bfK}{{\mbox{\boldmath $K$}} }

\newcommand{\bfM}{{\mbox{\boldmath $M$}} }

\newcommand{\bfP}{{\mbox{\boldmath $P$}} }

\newcommand{\bfQ}{{\mbox{\boldmath $Q$}} }

\newcommand{\bfF}{{\mbox{\boldmath $F$}} }

\newcommand{\bfI}{{\mbox{\boldmath $I$}} }

\newcommand{\bfJ}{{\mbox{\boldmath $J$}} }

\newcommand{\bfN}{{\mbox{\boldmath $N$}} }

\newcommand{\bfh}{{\mbox{\boldmath $h$}} }

\newcommand{\bfomega}{{\mbox{\boldmath $\omega$}} }

\def\XXint#1#2#3{{\setbox0=\hbox{$#1{#2#3}{\int}$ }
\vcenter{\hbox{$#2#3$ }}\kern-.6\wd0}}

\title{On the Motion of a Body with a Cavity Filled\\ with  Compressible Fluid}
\author{G. P. Galdi\, V. M\'acha\,  \v S. Ne\v casov\' a}
\input pomocne.in
\def\pat{\partial_t}

\begin{document}
\maketitle

\begin{abstract} We study the motion of the system, $\mathcal S$, constituted by a rigid body, $\mathcal B$, containing in its interior a viscous compressible fluid, and moving in absence of external forces. Our main objective is to characterize the long time behavior of the coupled system body-fluid. Under suitable assumptions  on the ``mass distribution" of $\mathcal S$, and for sufficiently ``small'' Mach number and initial data, we show that every corresponding motion (in a suitable regularity class) must tend to a steady state where the fluid is at rest with respect to $\mathscr B$. Moreover, $\mathcal S$, as a whole, performs a uniform rotation around an axis parallel to the (constant) angular momentum of $\mathcal S$, and passing through its center of mass.  
\end{abstract}

\section{Introduction}
The study of the motion of a rigid body with an interior cavity entirely filled with a liquid has a long-lasting history, tracing back to the pioneering contributions of Stokes \cite{Stokes}, Zhukovsky \cite{Zhu}, Poincar\'e \cite{Poin}, and, later, Sobolev \cite{Sob}.  

The principal and interesting characteristic arising in the dynamics of these coupled systems is that the presence of the liquid can substantially change the motion of the body and  may often produce  a significant stabilizing effect at  ``sufficiently large" times. 

It was  only in the late fifties, however, that this research area has become the object of a coherent and methodical investigation --mostly by Russian  mathematicians-- probably also due to the central role that it plays in rocket and space engineering, among the many other areas of concrete applications. The list of main contributors is too long to be included here, and we refer the interested reader to the classical monographs \cite{MoRu}, \cite{KKN}, \cite{Ch}, the latest book  \cite{Ch1} and the vast bibliography therein cited. 
\par
It must be emphasized that all the above contributions
are rarely of rigorous nature since 
they almost constantly rely on  simplified equations --some times, even ordinary differential
equations-- and/or special shapes of the body and cavity.
\par
A few years ago, a rigorous and systematic  mathematical analysis has eventually begun, under the assumption that the liquid filling  the cavity is  incompressible and governed by the Navier-Stokes equations   \cite{ST2,DGMZ,Ga,GMZ,GM,GMM,GMM1,GM1,MPS}.\footnote{ See also \cite{L1,L2,L3} for stability issues in the case of an inviscid liquid.}  Besides the basic study of well-posedness of the relevant initial-boundary value problem, all these works mainly focus on the fundamental question of the ``ultimate dynamics" of the coupled system, and the characterization of terminal states. \par In particular, in \cite{DGMZ} and \cite{Ga} the authors address the case where the system moves freely in absence of external forces and prove, among other things, that terminal states must be uniform rotations where the system, as a whole rigid body, spins with constant angular velocity around one of the central axes of inertia. This result --that is quite at odds with the classical one in absence of liquid-- holds in a very large class of motions (weak solutions) and whatever the shape of the body, of the cavity and the value of the physical parameters of the system, thus giving a rigorous proof of a famous conjecture formulated by {\sc N.Y. Zhoukovski} in 1885 \cite{Zhu}.    

The main purpose of this paper is to investigate whether and to what extent the above significant property continues to be valid if the fluid filling the cavity is assumed to be compressible (and isothermal). 

It must be remarked that, with the current knowledge of compressible Navier-Stokes fluids and  unlike  \cite{DGMZ,Ga}, it does not appear feasible   to  perform this type of analysis in the general class of weak solutions, whose initial data are only requested to have finite total energy. In fact, we leave this aspect of the problem as an outstanding open question. The reason for such a difficulty is easily seen. By its own nature, the problem calls for a detailed investigation of the $\Omega$-limit set associated to the generic weak solution. As recognized in \cite{DGMZ,Ga}, one of the crucial points  is to show that this set meets the basic property of  invariance. The latter, however, requires that weak solutions become as smooth for  large times as to ensure continuous dependence upon the data. In  \cite{DGMZ,Ga} this is shown to be the case  thanks to the ``smoothing property" possessed by weak solutions to the incompressible Navier-Stokes equations as times approaches infinity. 
In the analogous compressible problem, existence of a global weak solution for arbitrary initial
data possessing only finite energy  can be shown by a suitable modification of the arguments given in   \cite{LI4} and
 \cite{EF70} (see also \cite{FENO6}, \cite{NoSt}). However, with the current understanding of the theory of weak solutions, it is not possible to show, even with very regular initial data, that such a solution becomes smooth  for sufficiently  large times.  In this regards,  we recall that the propagation of singularities in weak solutions for a barotropic fluid  was studied in \cite{HofSan} in both two- and three-dimensional settings. There it is proved that, if the initial density is  nonnegative and essentially
bounded, initial energy is small, and initial velocities are in certain
fractional Sobolev spaces, then discontinuities in the density may persist for
all times; see \cite{Hoff} for the one-dimensional case.

In view of all the above, we  carry out our analysis within the class of ``strong'' solutions. The model we adopt is that of a barotropic gas where pressure and density are related by the formula $p=a\,\rho^\gamma$, $\gamma>1$, $a=\textrm{const.}>0$. Since we are interested in global existence, we shall need a restriction on the magnitude of the initial data in appropriate spaces; see Theorem 2. This allows us to show the existence of an $\Omega$-limit set possessing the desired properties of compactness, connectedness and invariance; see Lemma 7. As a consequence of the invariance, we then prove that this set is contained in the class of all possible  steady-states solutions. As shown in Proposition 1, each of these states is characterized by the fluid being at rest relative to the body, and  density distribution, $\rho_s$, and  constant angular velocity, $\omega_s$, satisfying a coupled and highly nonlinear system of equations; see (\ref{stationary2}). These equations translate the physical requirement that, since the fluid is at rest, the gradient of pressure must balance the centrifugal force imparted by the rotation which, in turn, occurs around an axis that depends on the density distribution itself. 
In fact, for a given admissible density distribution $\rho_s$ of the fluid, let $G=G(\rho_s)$ be the center of mass of the  coupled system body-fluid, $\mathcal S=\mathcal S(\rho_s)$, and let $\bfI=\bfI(\rho_s)$ be the inertia tensor of $\mathcal S$ relative to $G$ ({\em central inertia tensor}). Then,  the generic steady-state solution describes a motion where the coupled system $\mathcal S$  uniformly rotates, as a whole rigid body and with angular velocity $\omega_s=\omega_s(\rho_s)$, around an axis parallel to one of the eigenvectors of $\bfI$ and passing through $G$. 
\par
The final and critical property that remains to be investigated is under which  conditions  the $\Omega$-limit set corresponding to a given solution reduces to a singleton. As is well known, this property is by no means guaranteed, even for the ``classical"\footnote{ Namely, the body is kept at rest at all times, and the fluid moves in a bounded domain.} initial-boundary value problem for a barotropic fluid subject to a {\em given} potential-like force \cite{BdV,FePe,FP1,Er}.  As a matter of fact, in our case the situation appears to be more complicated. In fact, as mentioned in the previous paragraph, in any possible terminal state the gradient of pressure has to balance the potential-like centrifugal force where the angular velocity is itself an unknown function of the density. The problem of uniqueness  is addressed in Section 6 by means of the implicit function theorem, and the relevant findings are collected in Theorems 3 and 4. The results depend on the magnitude of the total angular momentum, $M_0$, of the coupled system $\mathcal S$. Since the motion occurs in absence of external forces, this quantity is preserved at all times; see (\ref{M_0_0}). If $M_0=0$ (a rather non-generic situation) then the $\Omega$-limit set reduces to the single state where $\mathcal S$ is at rest and the density of the gas is a constant (compatible with conservation of mass). If $M_0\neq 0$, the situation appears to be much more complicated. Nevertheless, we are able to give an answer  under suitable assumption on the ``mass distribution" of the coupled system, and by requiring the material constant $a$ to be ``sufficiently large" (``small'' Mach number).    
More precisely, let $\bar{\rho}\equiv m_{\mathcal F}/{\mathcal V}$, with $m_{\mathcal F}$ mass of the gas and $\mathcal V$ volume of the cavity. Then, $\bfI(\bar{\rho})$ is the central inertia tensor of the system body-fluid when the fluid has a constant density $\bar{\rho}$. Notice that $\bfI(\bar{\rho})$ is completely characterized by  known physical parameters. Then, if all three eigenvalues of $\bfI(\bar{\rho})$ are distinct, and $ a>a_0$ for a sufficiently large $a_0>0$, the $\Omega$-limit corresponding to any strong solution reduces to a single point or, equivalently, the generic strong solution to the initial-boundary value problem tends to a uniquely determined steady-state.    
\par
The paper is organized as follows. In Section \ref{FR}
we formulate the problem in both inertial and body-fixed frames and provide some preliminary considerations of general nature. The steady-state problem is addressed in Section \ref{steady.state}. In particular, it is shown that, in the very general class of weak solutions, the only steady-states are those where the fluid is at rest relative to the body.
Sections \ref{section.local.sol} and \ref{section.global.sol} are
devoted to the existence of  strong solutions. We follow the
approach developed by Matsumura and Nishida (\cite{MatNis},
\cite{MatNis2}) and  A. Valli \cite{valli}. The proof of the
existence of local-in-time strong solutions is performed, basically, along the same lines introduced in \cite{MatNis,valli}, with only a few differences in details . Thus, in Section
\ref{section.local.sol} we focus only on these differences.
However, the situation with global-in-time solutions is rather dissimilar from that addressed in the cited work,  due
to the presence of non-homogeneous boundary conditions. In \cite{valli}, where 
homogeneous Dirichlet boundary conditions are adopted, the author shows that 
in absence of external forces, as time $t\to\infty$ the velocity of the fluid tends to
zero and density tends to a constant function. This cannot be true
in our case, in general, since the magnitude of total angular momentum of the coupled system body-fluid
must be conserved (see (\ref{M_0})). Consequently, the velocity of the whole system cannot be zero in the limit $t\to\infty$ and also, contrary to the incompressible case, we cannot expect a constant density as a stationary solution.
In Section \ref{section.global.sol} we provide a proof of existence of global-in-time smooth solution whose velocity does not necessarily tend to zero.

The concluding section contains the main contribution of our paper,  devoted to long-time-behavior of strong solutions. We provide a number of properties  of $\Omega$-limit set, and show, in particular, that it must be  a subset of the class of steady-state solutions. Finally,  as  discussed earlier on in this section, we provide sufficient conditions ensuring that this set consists of one point only; see Theorems 3 and \ref{final.claim}.
\renewcommand{\theequation}{\arabic{section}.\arabic{equation}}
\section{Formulation of the Problem and Governing Equations}\label{FR}
Let $\B$ be a rigid body with an interior cavity, $\Ce$, entirely filled with a viscous,  compressible fluid, $\F$, and let $\S$ indicate the coupled system body-fluid (i.e. $\S:=\B\cup \F$). We are interested in the case where $\S$ moves in absence of external forces (inertial motions). This implies, in particular,  that its center of mass, $G$,  performs a uniform, rectilinear motion with respect to an inertial frame $\mathscr I$. Denoting by $\mathscr F$ the frame with the origin at $G$ and axes constantly parallel to those of $\mathscr I$, we have that $\mathscr F$ is inertial as well and that $G$ is at rest in $\mathscr F$. The  governing equations of $\S$ with respect to $\mathscr F$ then take  the following form
\begin{equation}\label{first.sys}
\begin{array}{cc}\medskip
\left.\begin{array}{rr}\medskip
\pat (r \bfw) + \diver (r \bfw \otimes \bfw) - \diver T(\bfw, p(r)) =
0\\
\pat r + \diver(r \bfw) =
0\end{array}\right\}\ \ \mbox{$(y,t)\in\cup_{t>0}\,\Ce(t)\times \{t\}$}\,, \\ \medskip
\bfw = 
\bfom(t)\times (\bfy-\bfy_{C}) + \bfeta(t)\,,\  \ \ \mbox{$(y,t)\in\cup_{t>0}\,\partial\Ce(t)\times \{t\}$}\,,\\ \medskip
\left.\begin{array}{rr}\medskip
\frac{d}{dt} (\bfJ_C \cdot\bfom) = -\int_{\partial \Ce(t)} (\bfy-\bfy_C)\times T(\bfw, p(r))\cdot \bfN\\
m_\B\,\bfeta(t) = -\int_{\mathcal C(t)} r\,\bfw
\end{array}\right\}\ \ t\in(0,\infty)\,,
\end{array}
\end{equation}
where $r$, $p$ and $\bfw$ are fluid density, pressure and velocity fields, $\bfom$ is the  angular velocity of $\B$, and $\bfeta$ the velocity of its center of mass $C$. Moreover, $\bfy_C$ denotes the vector position of $C$, while $m_\B$,  $\bfJ_C$ are mass and  inertia tensor  with respect to $C$ of $\B$, respectively, and $\bfN$ unit outer normal on $\partial\mathcal C$. Also, 
\begin{equation}\label{T}
T(\bfw,p) = S(\nabla \bfw) - {\bf1} p
\end{equation} 
(${\bf 1}$=\,unit tensor) is the Cauchy stress tensor with $S$ defined by 
\begin{equation}\label{S}S(\nabla \bfw) = \mu\D(\bfw) + (\lambda -\mbox{$\frac{2}{3}$} \mu) \bf1 \diver \bfw,\,,\ \ \D(\bfw):=\nabla \bfw+ (\nabla\bfw)^\top\,,\ \ \mu>{\rm 0}\,, \ \lambda\geq {\rm 0}\,.\end{equation} 
As for the dependence of $p$ on $r$, we shall consider the isentropic case 
\begin{equation}\label{pressure.rule}
p(r) = a\,r^\gamma,
\end{equation}
where $\gamma$ (specified later) and $a$ are positive material constants.

Equations (\ref{first.sys})$_{1,2,3}$ represent conservations of linear momentum and mass for $\mathcal F$, along with adherence of the fluid at the boundary of $\Ce$, whereas (\ref{first.sys})$_{4}$ is the balance of angular momentum of $\B$. Finally,   (\ref{first.sys})$_5$ translates the fact that the center of mass $G$ 
of $\S$ is at rest in $\mathscr F$.

As customary in this type of problems \cite{G2,KNN,NTT,MN,MN1}, it is convenient to rewrite the system (\ref{first.sys}) in a frame, $\mathscr R$, attached to $\B$, so that the domain occupied by the fluid becomes time-independent. To this end, let $\bfQ=\bfQ(t)$, $t\ge 0$, be the family of proper orthogonal transformations
defined by the equations
$$ \mbox{$\frac d{dt}$} \bfQ(t) = \mathbb S(\bfom)\cdot\bfQ(t)\,,\ \ \bfQ(0)={\bf 1}\,, 
$$ 
where 
\begin{equation}\label{S1}
\mathbb S (\bfa):=\left[\begin{array}{ccc}\smallskip 
0 & -a_3 & a_2\\ \smallskip
a_3 & 0 &-a_1\\
-a_2 &a_1& 0\end{array}\right]\,.
\end{equation}
By choosing $C$ as the origin of $\mathscr R$ we perform the following change of coordinates
$$
\bfx = \bfQ^\top\cdot(\bfy-\bfy_C)
$$
and define accordingly the transformed quantities
\begin{equation}
\begin{split}\label{transformation}
\rho(t,x) &= r(t,\bfQ^\top\cdot \bfx+\bfy_C)\,,\ \ \bfu(t,x) = \bfQ^\top \cdot\bfw(t,\bfQ^\top\cdot \bfx+\bfy_C)\\
\bfomega(t) &= \bfQ^\top\cdot\bfom(t)\,,\ \ \bfxi(t)=\bfQ^\top\cdot\bfeta(t)\,,\ \ \bfI_C = \bfQ\cdot\bfJ_C(t)\cdot \bfQ^\top\,,\ \ \bfn = \bfQ^\top\cdot\bfN(t)\,. 
\end{split}
\end{equation}
As a result, one shows that (\ref{first.sys}) becomes \cite{G2,KNN}
\begin{equation}\label{first.sys.1}
\begin{array}{cc}\medskip
\left.\begin{array}{rr}\medskip
\pat (\rho \bfu) + \diver (\rho\, \bfv \otimes \bfu) +\rho\,\bfomega\times\bfu+ \nabla p(\rho) =\diver S(\nabla\bfu)
\\
\pat \rho + \diver(\rho \bfv) =
0\end{array}\right\}\ \ \mbox{in \,$\Ce\times (0,\infty)$}\,, \\ \medskip
\bfu = 
\bfomega(t)\times \bfx + \bfxi(t)\,,\  \ \ \mbox{on\, $\partial\Ce\times (0,\infty)$}\,,\\ \medskip
\left.\begin{array}{rr}\medskip
\bfI_C\cdot\frac{d}{dt}\bfomega+\bfomega\times (\bfI_C\cdot\bfomega) = -\int_{\partial \Ce} \bfx\times T(\bfu, p(\rho))\cdot \bfn\\
m_\B\,\bfxi(t) = -\int_{\mathcal C} \rho\,\bfu
\end{array}\right\}\ \ t\in(0,\infty)\,,
\end{array}
\end{equation}
where
\begin{equation}\label{rel.vel}
\bfv:=\bfu-\bfomega\times\bfx-\bfxi\,,
\end{equation}
stands for the relative velocity field of the fluid with respect to the body. In this regard, we point out the obvious identity
$$
S(\nabla\bfu)=S(\nabla\bfv)\,,
$$
which will be often used throughout this paper without explicit mention.
\par
For future reference, we also observe that, integrating both sides of (\ref{first.sys.1})$_1$ over $\Ce$ --also by parts where needed-- and taking into account (\ref{T}) and the fact that $\bfv|_{\partial \Ce}={\bf 0}$, we easily get
$$
\frac d{dt}\int_{\Ce}\rho\,\bfu+\bfomega\times\int_{\Ce}\rho\,\bfu=\int_{\partial\Ce}T(\bfu,p(\rho))\cdot\bfn\,.
$$ 
The latter, in combination with (\ref{first.sys.1})$_5$ then furnishes
\begin{equation}\label{RM}
m_{\B}\frac {d\bfxi}{dt}+m_\B\bfomega\times\bfxi=-\int_{\partial\Ce}T(\bfu,p(\rho))\cdot\bfn\,.
\end{equation}
\par
Proceeding as in \cite[Appendix]{DGMZ}, one can further prove that (\ref{first.sys.1}) can be put in the following equivalent form that will turn out to be useful for our purposes
\begin{equation}\label{first.sys.2}
\begin{array}{cc}\medskip
\left.\begin{array}{rr}\medskip
\pat (\rho \bfu) + \diver (\rho\, \bfv \otimes \bfu) +\rho\,\bfomega\times\bfu+ \nabla p(\rho) =\diver S(\nabla\bfu)
\\
\pat \rho + \diver(\rho \bfv) =
0\end{array}\right\}\ \ \mbox{in \,$\Ce\times (0,\infty)$}\,, \\ \medskip
\bfu = 
\bfomega(t)\times \bfx + \bfxi(t)\,,\  \ \ \mbox{on\, $\partial\Ce\times (0,\infty)$}\,,\\ \medskip
\left.\begin{array}{rr}\medskip
\frac{d}{dt}\bfM+\bfomega\times \bfM ={\bf 0},\ \ \bfM:=\bfI_C\cdot\bfomega+\int_{\Ce}\rho\,\bfx\times\bfu\\
m_\B\,\bfxi(t) =-\int_{\Ce}\rho\,\bfu\,
\end{array}\right\}\ \ t\in(0,\infty)\,.
\end{array}
\end{equation}
Notice that (\ref{first.sys.2})$_4$ states conservation of the angular momentum of the coupled system $\S$ with respect to center of mass of $\B$.
Finally, we endow equations (\ref{first.sys.1}) (or, equivalently, (\ref{first.sys.2}))  with the following initial conditions
\begin{equation}\label{initcon}
\rho(0) = \rho_0,\qquad \rho(0)\bfu(0) = {\bfq}.
\end{equation}
\par
Notice also that from (\ref{first.sys.1})$_2$ and (\ref{initcon}) we formally deduce the following  global form of conservation of mass
\begin{equation}\label{conmass}  
\int_{\Ce}\rho(x,t)\,{\rm d}x=\int_{\Ce}\rho_0(x)\,{\rm d}x=m_\F\,,\ \ t\ge 0\,,
\end{equation}
where $m_\F$ is the total mass of the fluid.\par
The problem we shall address in the following sections regards the asymptotic behavior in time of solutions $(\rho,\bfu,\bfomega,\bfxi)$ to (\ref{first.sys.1}) in a suitable function class, $\mathscr C$ (say). More precisely, we shall show that, provided the size of the initial data is suitably restricted and under certain assumptions on the distribution of masses, any such solution in $\mathscr C$ will converge to a steady-state where the fluid is at relative rest ($\bfv\equiv{\bf 0})$  and $\S$  performs a uniform rotation as a single rigid body.

We end this section by recalling some basic notation  we will use throughout.
As customary, $L^p$ (resp.  $W^{k,p}$) denote the classical Lebesgue (resp. Sobolev) space,  with corresponding norm  $\|\cdot\|_p$ (resp. $\|\cdot\|_{k,p}$). Norms in Bochner space $L^p([0,T],L^q)$ (resp. $L^p([0,T], W^{k,q})$) are denoted by $\|\cdot\|_{L^p(L^q)}$ (resp. $\|\cdot\|_{L^p(W^{k,q})}$). Supporting sets in all the above spaces will be usually omitted. We also need the following  class of functions defined on the whole of $\S$ and that reduce to ``rigid motions" on $\B$:
\begin{equation*}
W^{k,p}_R(\S) = \{{\bfphi}:\S\mapsto \mathbb R^3, {\bfphi}\in W^{k,p}(\S), {\bfphi}|_{\B} = {\bfell}\times x  + {\bfk}\mbox{ for some }{\bfell}, {\bfk}\in \mathbb R^3\}\,.
\end{equation*}
We set $L^p_{R}(\S):= W^{0,p}_R(\S)$.
If $\bfphi\in W^{k,p}_R(\S)$ we shall typically write
$\bfell_{\mbox{\footnotesize{$\bfphi$}}}$ and $\bfk_{\mbox{\footnotesize{$\bfphi$}}}$ to emphasize that the characteristic vectors $\bfell$ and $\bfk$ are associated to $\bfphi$. For future reference,  we observe that, due to the  properties of the tensor field $S$ defined in (\ref{S}), for any $\bfpsi,\bfphi\in W^{1,2}_R(\S)$ we have   
$$
S(\nabla \bfpsi):\nabla \bfphi = S(\nabla \bfpsi):\nabla(\bfphi - \bfell_{\mbox{\footnotesize{$\bfphi$}}}\times \bfx),
$$
and also
$$
S(\nabla \bfphi):\nabla \bfphi = S(\nabla \bfphi- \bfell_{\mbox{\footnotesize{$\bfphi$}}}\times \bfx)):\nabla (\bfphi- \bfell_{\mbox{\footnotesize{$\bfphi$}}}\times \bfx).
$$

We finally remark that, in what follows, the cavity $\mathcal C$ is assumed to be a bounded domain of $\mathbb R^3$ possessing (at least) Lipschitz regularity. 
\setcounter{equation}{0}
\section{Steady-State Solutions}
\label{steady.state}
Objective of this section is to derive certain  basic properties of steady-state solutions to  (\ref{first.sys.2}). Sufficient conditions for their existence will be postponed till the last section.
\par We begin to observe that 
from (\ref{first.sys.2}) we have that steady-state solutions must satisfy the following set of equations
\begin{equation}\label{stationary}
\begin{array}{cc}\medskip
\left.\begin{array}{rr}\medskip\diver(\rho\bfv\otimes\bfu) + \rho\omega\times \bfu + \nabla p(\rho)=\diver S(\nabla \bfu)\\
\diver{\rho \bfv}=0\end{array}\right\}\,\  \mbox{on}\ \mathcal C\\ \medskip
\bfu=\bfomega\times \bfx + \bfxi\ \mbox{on}\ \partial \mathcal C\\ \medskip
m_\B\bfxi = -\int_{\Ce}\rho\bfu\\
\bfomega\times \bfM = 0\,,\ \ \bfM:=\bfI_C\cdot\bfomega+\int_{\Ce}\rho\,\bfx\times\bfu.
\end{array}
\end{equation}

If we formally dot-multiply (\ref{stationary})$_1$ by $\bfphi\in C_0^\infty(\mathcal C)$ and integrate by parts over $\mathcal C$, we get  
\begin{equation}\label{wfm}
\int_{\Ce}\big[-(\rho\,\bfv\otimes\bfu):\nabla\bfphi+\rho\,\bfomega\times\bfu\cdot\bfphi+p(\rho)\diver\bfphi\big]=-\int_{\Ce}S(\nabla\bfu):\nabla\bfphi\,,\ \ \mbox{all $\bfphi\in C_0^\infty(\mathcal C)$}\,.
\end{equation}
Likewise, from (\ref{stationary})$_2$ we show
\begin{equation}\label{wfcm}
\int_{\Ce}\rho\,\bfv\cdot\nabla\psi=0\,,\ \ \mbox{all $\psi\in C_0^\infty(\mathcal C)$}\,.
\end{equation}
Moreover, if 
\begin{equation}\label{assumpt1}
b\in C^1(\mathbb R_+),\, b(0)=0, b'(r) \geq C_b,
 \end{equation}
 with $C_b$ a constant that may depend on $b$,  again from (\ref{stationary})$_2$ we deduce the ``renormalized" continuity equation
\begin{equation}\label{rn}
\diver(b(\rho)\bfv)+\big(\rho\,b^\prime(\rho)-b(\rho)\big)\diver\bfv=0\ \ \mbox{in $\Ce$}\,,\ \ \mbox{all $b\in C^1(\mathbb R_+)$}\,.
\end{equation}

We then say that the quadruple $(\rho,\bfu,\bfomega,\bfxi)$ is a {\em renormalized weak solution} to (\ref{stationary}) if, for some $\gamma>1$, (i) $(\rho,\bfu) \in L^\gamma(\mathcal C)\times W^{1,\gamma^*}_R(\mathcal S)$, $\bfomega=\bfomega_{\mbox{\footnotesize{$\bfu$}}}$, $\bfxi=\bfxi_{\mbox{\footnotesize{$\bfu$}}}$,\footnote{Here $\gamma^* = \frac{3\gamma}{3+\gamma}$ for $\gamma > 6$, $\gamma^* = \frac{9\gamma}{5\gamma - 3}$ for $\frac32<\gamma\leq6$, $\gamma^* = 3+$ for $\gamma = \frac32$, and $\gamma^* = \frac{\gamma}{\gamma-1}$ for $\gamma<\frac32$.}\,(ii) $(\rho,\bfu,\bfomega,\bfxi)$ satisfies (\ref{stationary})$_{4,5}$, (\ref{wfm}), (\ref{wfcm}), and, in addition, (\ref{rn}) in the sense of distributions in the whole of $\mathbb R^3$, with $\rho$ and $\bfv$ prolonged by 0 outside $\Ce$ with $b$ satisfying (\ref{assumpt1}) and the following assumptions

\begin{equation}
\label{assumpt2}
\begin{array}{l}
|b'(t)|\leq c t^{-\lambda _0}, t \in (0,1], \lambda _0<1,\\
|b'(t)|\leq ct^{\lambda _1}, t\geq 1,c>0, -1<\lambda_1<\infty.
\end{array}
\end{equation}
The next lemma shows that weak solutions may occur only if the fluid is at rest relative to $\B$, namely, the coupled system $\S$ moves, as a whole, by rigid motion. 
\Bl\label{v0}
Let  $(\rho,\bfu,\bfomega,\bfxi)$ be a renormalized weak solution to (\ref{stationary}). Then $\bfv\equiv 0$.
\end{lemm}
\begin{proof} Clearly, $\bfv\in W_0^{1,\gamma^*}(\Ce)$. Since $C_0^\infty(\Ce)$ is dense in $W_0^{1,\gamma^*}(\Ce)$ we can replace $\bfv$ for $\bfphi$ in (\ref{wfm}) to get
\begin{equation}\label{wfm1}
\int_{\Ce}\big[-(\rho\,\bfv\otimes\bfu):\nabla\bfv+\rho\,\bfomega\times\bfu\cdot\bfv+p(\rho)\diver\bfv\big]=-\int_{\Ce}S(\nabla\bfu):\nabla\bfv\,.
\end{equation}
We now observe that
\begin{multline}\label{p1}
\int_{\mathcal C}\diver(\rho \bfv \otimes \bfu)\cdot \bfv = \int_\Ce \diver(\rho\,\bfv\otimes \bfu)\bfu - \int_\Ce \diver(\rho\bfv\otimes\bfu)(\omega\times x + \xi)\\ =  -\frac 12\int_{\mathcal C} \rho \bfv \nabla|\bfu|^2 - \int_{\mathcal C}\rho \omega \times \bfu \bfv = -\int_{\Ce} \rho\omega\times \bfu\bfv\,.
\end{multline}
Furthermore, choosing $b(\rho)=\rho\int_1^\rho\frac{p(s)}{s^2}ds$ we infer
$$
\rho\,b^\prime(\rho)-b(\rho)=p(\rho)\,,
$$
so that from (the distributional form of) (\ref{rn}) we show
\begin{equation}\label{p2}
\int_{\mathcal C}p(\rho)\diver\bfv = -\int_{\mathcal C} \diver\big(b(\rho)\bfv\big) =  0\,.
\end{equation}
Collecting  (\ref{wfm})--(\ref{p2}) we thus obtain $\int_{\mathcal C} S(\nabla \bfv):\nabla \bfv = 0$, which together with Korn's inequality yields the desired claim.
\QED\end{proof}
Setting
$$
\bar{\bfI}=\bar{\bfI}(\rho):=\int_{\Ce}\rho[{\bf 1}\,|\bfx|^2-\bfx\otimes\bfx]\,,
$$
from the previous lemma and (\ref{stationary}) we easily show that any weak solution to (\ref{stationary})  must be then of the form $(\rho_s,\bfv\equiv{\bf0},\bfomega_s,\bfxi_s)$, with $\rho_s$, $\bfomega_s$, and $\bfxi_s$ satisfying the following system of equations:
\begin{equation}\label{sasi}\begin{array}{ll}\medskip
\rho_s[\bfomega_s\times (\bfomega_s\times\bfx_s+\bfxi_s)]=-\nabla p(\rho_s)\\ \medskip
\,\bfxi_s=-{\displaystyle{\frac1{m_{\S}}}}\int_{\Ce}\rho_s\, \bfomega_s\times\bfx\\
\bfomega\times\big[(\bfI_C+\tilde{\bfI}_C)\cdot\bfomega+\int_{\Ce}\rho_s\,\bfx\times \bfxi_s\big]={\bf 0}\,,
\end{array}
\end{equation}
where $m_{\S}$ is the mass of $\S$ and $\tilde{\bfI}_C:=\bar{\bfI}(\rho_s)$
is (for $\rho_s>0$) the inertia tensor with respect to $C$ of the fluid in the steady-state configuration. \par
We notice that, if $\rho_s(x)>0$ in $\Ce$,  then by a simple boot-strap argument  from (\ref{sasi})$_1$ it follows that, in fact, $\rho_s\in C^\infty(\mathcal C)$.\par
Set
\begin{equation}\label{g}
\bfg=\bfg(\rho):=\int_{\Ce}\rho\,\bfx\,,
\end{equation}
and define 
$$
\bfI_g:=\frac1{m_{\S}}\left({\bf 1}|\bfg|^2-\bfg\otimes\bfg\right)\,.
$$
We have the following.
\Bl\label{stab} 
For any (sufficiently smooth) $\rho=\rho(x)>0$, the tensor
\begin{equation}\label{I}
\bfI=\bfI(\rho):= \bfI_C+\bar{\bfI}-\bfI_g\,.
\end{equation}
is symmetric and positive definite. Moreover,  denoting by $G=G(\rho)$ the center of mass of $\mathcal S$, $\bfI(\rho)$ coincides with the inertia tensor of $\mathcal S$ with respect to $G$.\El
\begin{proof} The symmetry property is obvious. Moreover, by using Lagrange vectorial identity, for any $\bfa\in\mathbb R^3$ we have
$$\begin{array}{rl}\medskip
m_{\S}\,\bfa\cdot\bfI_g\cdot\bfa&=|\bfa|^2|\bfg|^2-(\bfa\cdot\bfg)^2=(\bfa\times\bfg)^2=\big(\int_{\Ce}\rho\,\bfa\times\bfx\big)^2\\
&\leq \int_{\Ce}\rho\cdot\int_{\Ce}\rho(\bfa\times\bfx)^2=m_{\mathcal F}\bfa\cdot\bar{\bfI}(\rho)\cdot\bfa
\,,\end{array}
$$
where $m_{\mathcal F}:=\int_{\Ce}\rho$ is the mass of the fluid. Therefore,
$$
\bfa\cdot\bfI\cdot\bfa\ge\bfa\cdot\bfI_C\cdot\bfa+\big(1-\frac{m_{\mathcal F}}{m_{\S}}\big)\bfa\cdot\bar{\bfI}\cdot\bfa\,,
$$
which proves the desired property since $m_{\S}>m_{\mathcal F}$. We next observe that, denoting by $\bfx_G$ the vector position of $G$ in the frame $\mathscr R$, $\mathcal B$ the volume occupied by the body, and $\rho_{\mathcal B}$ its density, we have
$$
\bfx_G=\frac1{m_{\mathcal S}}\left(\int_{\mathcal B}\rho_{\mathcal B}\,\bfx+\int_{\mathcal C}\rho\,\bfx\right)=\frac1{m_{\mathcal S}}\,\bfg\,,
$$
since $
\int_{\mathcal B}\rho_{\mathcal B}\,\bfx=0 
$. The claimed property then follows from  classical Steiner's theorem on the variation of the inertia tensor with the pole, e.g., \cite[Section 3.5]{Scheck}.
\QED\end{proof}

Collecting all the above results we can now give the following characterization of the class of weak solutions to (\ref{stationary}). 
\Bp
\label{prop}
Let  $(\rho_s,\bfu_s,\bfomega_s,\bfxi_s)$ be a weak solution to (\ref{stationary}) with $\rho_s>0$. Then 
\begin{equation}\label{u}\bfu_s=\bfomega_s\times\bfx+\bfxi_s\,,\ \ x\in \S
\,,\end{equation} 
while $\rho_s$, $\bfomega_s$ and $\bfxi_s$ satisfy the following equations
\begin{equation}\label{stationary2}
\begin{array}{ll}\medskip
\rho_s^{\gamma-1}(x)  = \frac {\gamma - 1}{2a\gamma}\left(
|\bfomega_s\times \bfx|^2 -2 (\bfomega_s\times \bfxi_s)\cdot\bfx\right) + c
\,,\  x\in \Ce\,,\  \mbox{some } c\in \mathbb R,\\ \medskip
\bfomega_s\times (\bfI(\rho_s)\cdot \bfomega_s) = {\bf 0}\,,\\
m_\S \bfxi_s  = -\bfomega_s\times \bfg(\rho_s)\,.
\end{array}
\end{equation}
\Ep
\begin{proof} The property in (\ref{u}) follows from Lemma \ref{v0}. Moreover (\ref{stationary2})$_3$ is a consequence of (\ref{sasi})$_2$ and the definition (\ref{g}). In addition,
from (\ref{sasi})$_1$ we at once deduce that
\begin{equation*}
\rho_s \nabla\left(\frac 12|\bfomega_s\times \bfx|^2 - (\bfomega_s\times \bfxi_s)\cdot \bfx\right) = \nabla p(\rho_s)\mbox{ in }\Ce.
\end{equation*}
which furnishes (\ref{stationary2})$_1$. Finally, using (\ref{stationary2})$_{3}$ we get
$$
\int_{Ce}\rho_s\bfx\times\bfxi_s=\frac 1{m_\S}\,\bfg\times(\bfg\times \bfomega_s)=\frac 1{m_\S}\big[(\bfomega_s\cdot\bfg)\bfg-\bfomega_s\,|\bfg|^2\big]=-\bfI_g(\rho_s)\cdot\bfomega_s\,,
$$
which, in conjunction with (\ref{sasi})$_3$, produces (\ref{stationary2})$_2$. This concludes the proof.
\QED\end{proof}

\Br
As  mentioned at the beginning of this section, the existence  of  solutions to (\ref{stationary2}) (or, equivalently, weak solutions to (\ref{stationary})) will be addressed in the last section. More specifically, we shall show that  the $\Omega-$limit set associated to (strong) solutions to (\ref{first.sys.2}) is not empty and that its generic element is a steady-state solution; see Remark \ref{bu}.  
\Er
\setcounter{equation}{0}
\section{Local existence of strong solutions}\label{section.local.sol}Our next objective is to show that (\ref{first.sys.2})--(\ref{initcon}) is solvable in a class of ``strong" solutions, at least in some open interval of time. 
To this end, we begin to put (\ref{first.sys.2}) in a suitable weak form where the involved field variables are defined on the whole system $\S$. Let \begin{equation}\label{s}\overline \rho = \frac 1{|\Ce|} \int_\Ce \rho_0\,,\ \ \sigma := \rho - \overline \rho\,,\end{equation} and consider the following set of equations 
\begin{equation}
\begin{split}\label{local.existence.strong}
\pat \sigma  + \bfv\cdot \nabla \sigma  + \sigma \diver \bfv + \overline \rho \diver \bfv &= 0\mbox{ in }\Ce\times (0,\infty),\\
\int_{\Ce}\sigma(x,t)\,{\rm d}x&= 0\mbox{ in }(0,\infty)\\
\int_\S \rho \pat \bfu\cdot \bfphi + \int_\S \rho\bfv\cdot\nabla \bfu\cdot \bfphi + \int_\S \rho \bfomega_{\mbox{\footnotesize{$\bfu$}}} \times \bfu \cdot\bfphi - \int_\S p(\rho)\diver \bfphi + \int_\S S(\nabla \bfu):\nabla \bfphi &= 0,\\
\mbox{ for all }& \bfphi \in W^{1,2}_R(\S),\\
m_\B \bfxi_{\mbox{\footnotesize{$\bfu$}}} = -\int_{\Ce}\rho\bfu\,,\ \ t\in(0,\infty)\,.
\end{split}
\end{equation}
It is easy to see that every sufficiently smooth solution to (\ref{local.existence.strong}) is, in fact, a solution to (\ref{first.sys.1}), which, as shown earlier on, is equivalent to (\ref{first.sys.2}).  Actually, choosing at first $\bfphi\in C_0^\infty(\Ce)$,  integrating by parts (\ref{local.existence.strong})$_2$ as necessary  and taking into account (\ref{local.existence.strong})$_{1,3}$ we at once obtain that $\rho, \bfu$ solve (\ref{first.sys.1})$_{1,2,5}$. If we now take $\bfphi = {\bfell}\times \bfx$ and recall that $\S=\B\cup\Ce$ and that $\bfv\equiv{{\bf 0}}$ in $\mathcal C$, by a straightforward calculation  we show
$$
\int_{\Ce}\rho\, \big(\pat \bfu + \bfv\cdot\nabla \bfu + \bfomega_{\mbox{\footnotesize{$\bfu$}}} \times \bfu\big) \cdot\bfphi +\big[\frac{d}{dt} (\bfI_C\cdot\bfomega_{\mbox{\footnotesize{$\bfu$}}}) + \bfomega_{\mbox{\footnotesize{$\bfu$}}} \times (\bfI_C\cdot \bfomega_{\mbox{\footnotesize{$\bfu$}}})\big]\cdot\bfell=0\,.
$$
We next use (\ref{first.sys.1})$_{1,2}$ in the first integral to get
$$
\int_{Ce}\diver T(\bfu,p)\cdot\bfphi=-\big[\frac{d}{dt} (\bfI_C\cdot\bfomega_{\mbox{\footnotesize{$\bfu$}}}) + \bfomega_{\mbox{\footnotesize{$\bfu$}}} \times (\bfI_C\cdot \bfomega_{\mbox{\footnotesize{$\bfu$}}})\big]\cdot\bfell\,.
$$
The latter, in turn, after integration by parts delivers
$$
\big[\frac{d}{dt} (\bfI_C\cdot\bfomega_{\mbox{\footnotesize{$\bfu$}}}) + \bfomega_{\mbox{\footnotesize{$\bfu$}}} \times (\bfI_C\cdot \bfomega_{\mbox{\footnotesize{$\bfu$}}})\big]\cdot\bfell=-\bfell\cdot\int_{\partial\Ce}\bfx\times T(\bfu,p)\cdot\bfn\,,
$$
which, by the arbitrariness of $\bfell$ proves (\ref{first.sys.1})$_4$.
\smallskip\par
The major result in this section reads as follows. \Bt\label{local.strong.solution} Let $ \Ce$ be of class $C^3$,
$\bfu_0\in W^{1,2}_R(\S)$, $\bfu_0|_{\Ce}\in W^{2,2}(\Ce)$,
$\rho_0|_{\Ce}\in W^{2,2}(\Ce)$, $0< m \leq \rho_0 \leq M $ and
$\gamma >1$. Then there exists $T^*>0$,  $\bfu\in C([0,T^*],
W^{1,2}_R(\S))$, $\bfu|_{\Ce}\in L^2(0,T^*, W^{3,2})\cap C([0,T^*],
W^{2,2})$ with $\pat \bfu \in L^2(0,T^*, W^{1,2}_R(\S)) \cap
C([0,T^*], L^2_R(\S))$ and $\rho\in C([0,T^*], W^{2,2}(\Ce))$ with
$\pat \rho \in C([0,T*], W^{1,2}(\Ce))$, $\rho >0$ in $\Ce\times
[0,T^*]$ such that $(\bfu, \rho)$ is a solution to (\ref{s}),
(\ref{local.existence.strong}). Moreover, this solution is unique in its own 
class. \Et 

The proof of this theorem will be achieved by the  method employed in \cite{valli}. The main ingredients
are the Schauder fix point argument combined with
regularity results for the continuity equation and a suitable elliptic
problem. We split the proof  into three steps, each one described in the next three subsections.

\subsection{Continuity equation}
Consider the  problem
\begin{equation}
\begin{split}\label{continuity.equation.strong}
\pat \sigma + \tilde \bfv\cdot \nabla \sigma + \sigma \diver \tilde\bfv + \overline \rho \diver \tilde \bfv &= 0 \mbox{ in }\Ce\times(0,T)\\
\sigma(0) &= \sigma_0,
\end{split}
\end{equation}
where $\tilde \bfv$ and $\sigma_0$ are given functions. The following result holds
\Bl
\label{continuity.equation.lemma}
Let $\partial \Omega \in C^1$, $\tilde \bfv \in L^1(0,T;W^{3,2}(\Ce))$, $\tilde \bfv \cdot n = 0$ and $\sigma_0\in W^{2,2}(\Ce)$ with $\int_\Ce \sigma_0 = 0$. Then there exists a unique solution $\sigma$ to (\ref{continuity.equation.strong}) such that $\sigma \in C([0,T], W^{2,2}(\Ce))$, $\int_\Ce \sigma(t) = 0$ for each $t\in [0,T]$ and
\begin{equation*}
\|\sigma\|_{L^\infty(W^{2,2})}\leq c(\|\sigma_0\|_{W^{2,2}} + 1 ) \exp(\|\tilde \bfv\|_{L^1(W^{3,2})}).
\end{equation*}
Moreover, if $\tilde \bfv \in C([0,T], W^{2,2}(\Ce))$, then $\pat \sigma \in C([0,T], W^{1,2}(\Ce))$ and
\begin{equation*}
\|\pat \sigma\|_{L^\infty(W^{1,2})} \leq c \|\tilde \bfv\|_{L^\infty (W^{2,2})}(\|\sigma_0\|_{W^{2,2}} + 1 ) \exp(\|\tilde \bfv\|_{L^1(W^{3,2})}).
\end{equation*}
\El
\begin{proof}
See Lemma 2.3 in \cite{valli}.
\QED\end{proof}

\subsection{Balance of linear momentum}
In this subsection we are concerned with the problem
\begin{equation}\label{bolm.strong}
\begin{split}
\int_\S \tilde \rho\, \pat \bfw\cdot  \bfphi  + \int_\S S(\nabla \bfw): \nabla \bfphi &= [\bfF,\bfphi]\\
\bfw(0)&=\bfw_0\\
\bfxi_{\mbox{\footnotesize{$\bfu$}}} &= -\frac1{m_\S}\int_\Ce \tilde \rho \bfw.
\end{split}
\end{equation}
where $\bfphi \in W^{1,2}_R(\S)$, $\bfphi|_\B = \bfomega_{\mbox{\footnotesize{$\bfphi$}}}\times \bfx$, $\tilde \rho$ and $\bfF\in (W^{1,2}_R(\S))^*$ are given and $\bfu = \bfw + \bfxi_{\mbox{\footnotesize{$\bfu$}}} = \bfv + \bfomega_{\mbox{\footnotesize{$\bfu$}}}\times \bfx + \bfxi_{\mbox{\footnotesize{$\bfu$}}}$ for some (unknown) $\bfomega_{\mbox{\footnotesize{$\bfu$}}}$ and $\bfxi_{\mbox{\footnotesize{$\bfu$}}}$. Here $[.,.]$ denotes a duality pairing between $W^{1,2}_R(\S)$ and $(W^{1,2}_R(\S))^*$

\Bl\label{bolm.strong.lemma}
Let $\partial \Ce \in C^2$, $\tilde \rho \in L^\infty(( \S\times(0,T))$, $0\leq \frac m 2\leq \tilde \rho \leq 2M$ a.e. in $ \Ce\times [0,T]$, $\bfF\in L^2(0,T, L^2_R(\S))$, 
and $\bfu_0\in W^{1,2}_R(\S)$. Then the following properties hold.
\begin{itemize}
\item[{\rm (i)}]
There exists a unique solution $\bfu$ to (\ref{bolm.strong}) such that
$$
\bfu \in C(0,T,W^{1,2}_R(\S)), \bfu|_{\Ce} \in L^2(0,T, W^{2,2}(\Ce)), \ \pat \bfu \in L^2(0,T,L^2_R(\S)).
$$
\item[{\rm (ii)}] Suppose, in addition, $\partial \Ce \in C^3$, $0<m\leq \tilde \rho(0,x)\leq M$ a.e. in $\Ce$, $\nabla \left(\tilde \rho|_{\Ce}\right) \in L^4(0,T, L^6(\Ce))$, $\pat \tilde \rho \in L^2(0,T,L^3(\Ce))$, and  that $\bfF$, $\pat\bfF$ satisfy the following further assumptions
$$
[\bfF,\bfphi] = \int_\S \bfF_1\cdot \bfphi + \int_\S \bfF_2 : \nabla \bfphi,
$$
where $\bfF_1\in L^2(0,T, W^{1,2}_R(\S))\cap L^\infty(0,T,L^2)$ and $\bfF_2 \in L^2(0,T,W^{2,2}(\S))\cap L^\infty(0,T, W^{1,2}(\S))$, and
$$[\pat \bfF,\bfphi] = \int_\S \pat \bfF_3 \cdot\bfphi + \int_\S \pat \bfF_4 :\nabla \bfphi,$$
where $\pat \bfF_3\in L^2(0,T,L^{\frac 65}_R(\S))$ and $\pat \bfF_4 \in L^2(0,T,L^2_R(\S))$. Finally, assume $\bfu_0 \in W^{1,2}_R(\S)$ with $\bfu_0|_{\Ce}\in W^{2,2}(\Ce)$. Then
$$
\bfu|_{\Ce} \in L^2(0,T, W^{3,2}(\Ce))\cap C(0,T,W^{2,2}(\Ce)), \ \ \pat \bfu \in L^2(0,T, W^{1,2}_R(\S))\cap L^\infty(0,T, L^2_R(\S))
$$
and
\begin{multline*}
\|\bfu\|^2_{L^\infty(L^{2})} + \|\bfv\|_{L^\infty(W^{2,2})}^2 + \|\bfu|_{\Ce}\|_{L^2(W^{3,2})}^2 + \|\pat \bfu\|_{L^\infty(L^2)}^2 + \|\pat \bfu\|_{L^2(W^{1,2})}^2\\
\leq c(\|F_1\|_{L^2(W^{1,2})}^2 + \|\bfF_1\|_{L^\infty(L^2)}^2 + \|\bfF_2\|_{L^\infty(W^{1,2})}^2 + \|\bfF_2\|^2_{L^2(W^{2,2})})\\
 + c\left(\|\pat \bfF_3\|_{L^2(L^\frac 65)}^2 + \|\pat \bfF_4\|_{L^2(L^2)} + \|\bfu_0\|_{W^{1,2}}^2 + \|\bfF(0)\|_{L^2}^2\right)\\
\cdot (1+ \|\nabla \tilde \rho\|_{L^4(L^6)}^4 + \|\pat \rho\|_{L^2(L^3)}^2)\exp (c \|\pat \tilde\rho\|_{L^2(L^3)}^2)
\end{multline*}
\end{itemize}
\El
\begin{proof}
{\em Property {\rm (i)}}. 
We begin to derive some a priori estimates. We choose $\bfphi  = \pat \bfw$ in (\ref{bolm.strong}), and after some calculations we deduce
$$
\frac m 2\int_\S |\pat \bfw|^2  + c \pat \int_\Ce\left( |{\mathcal D} (\bfv) |^2 +  |\diver \bfv|^2\right) \leq \|F\|_{2}^2.
$$
By integrating this identity over $(0,t)$,  we derive  $\pat \bfw \in L^2(0,T, L^2(\S))$, $ \nabla \bfv \in L^\infty(0,T,L^2)$. Consequently, we get $\bfomega_{\mbox{\footnotesize{$\bfu$}}} \in L^\infty(0,T)\cap W^{1,2}(0,T)$, $\bfxi_{\mbox{\footnotesize{$\bfu$}}}\in L^\infty(0,T)$ and $\nabla \bfu \in L^\infty(0,T, L^2)$. \\
Also, from (\ref{bolm.strong}) we get 
$$
\int_\S S(\nabla \bfv): \nabla \bfphi = \int_\S \tilde \rho (-\pat \bfw)\cdot \bfphi + [\bfF,\bfphi]
$$
for every $\bfphi \in W^{1,2}_0(\Ce)$. Recalling that $\bfv|_{\partial \Ce} = 0$, thanks to classical regularity results for elliptic problems we get for a.a. $t$
$$
\|\nabla^2 \bfv\|_{2} \leq c \left(\|\tilde\rho\, \pat \bfw\|_{2} + \|\bfF_1\|_{2} + \|\bfF_2\|_{W^{1,2}}\right),
$$
and, consequently,
$$
\|\nabla^2 \bfv\|_{L^2(L^2)} \leq c \left(\|\tilde\rho \pat \bfw\|_{L^2(L^2)} + \|\bfF_1\|_{L^2(L^2)} + \|\bfF_2\|_{L^2(W^{1,2})}\right).
$$
With the above estimates in hand, the proof of  existence of unique solution to (\ref{bolm.strong}) is quite standard. In fact, for $\tilde \rho$ constant it can be achieved by means of the classical Galerkin approximation (see e.g. \cite[Chapter 7]{evans}). In the general case of $\tilde \rho$ not constant,  we may employ exactly the same approach as in  \cite[Lemma 2.1]{valli}. Note, that due to the assumed compatibility condition, the problem  is not over-determined. Thus the proof of Property (i) can be considered complete.
\smallskip\\
{\em Propery {\rm (ii)}}.
We choose in (\ref{bolm.strong})  $\bfphi = -\pat \pat \bfw$, and  integrate over $(0,t)$. In what follows we use the notation $\bfW := \pat \bfw$ and $\bfV := \pat \bfv$. After some straightforward calculations we get
\begin{multline}\label{derivace.U}
\frac 12\int_\S \tilde \rho |\bfW|^2 (t) + \int_0^t \int_\S S(\nabla \bfW):\nabla \bfW =\\
 \frac 12\int_\S \rho |\bfW|^2(0) + \int_0^t \int_\S \pat \bfF_3 \cdot\bfW  + \int_0^t \int_\S \pat \bfF_4 : \nabla \bfW + \frac 12 \int_0^t \int_\S\pat \tilde \rho |\bfW|^2
\end{multline}
Since all norms in finite-dimensional space are equivalent, we have 
\begin{equation}\label{10}
\|\mbox{$\frac{d}{dt}$}\bfomega_{\mbox{\footnotesize{$\bfu$}}}\times \bfx\|^2_{\frac 83}\sim |\mbox{$\frac{d}{dt}$}\bfomega_{\mbox{\footnotesize{$\bfu$}}}|^2\sim \int_\B \tilde\rho|\mbox{$\frac{ d}{dt}$}\bfomega_{\mbox{\footnotesize{$\bfu$}}} \times \bfx|^2 \leq c\int_\S \tilde\rho|\bfW|^2.\end{equation}
Consequently,
\begin{multline}
\label{odhad.nabla.U}
\int_0^t \int_\S |\pat \tilde\rho||\bfW|^2 \leq c\left(\int_0^t \int_\S |\pat \tilde\rho||\bfV|^2 + \int_0^t \int |\pat \tilde \rho||\mbox{$\frac d{dt}$}\bfomega_{\mbox{\footnotesize{$\bfu$}}}\times \bfx|^2\right) \\
\leq c\left(\int_0^t \|\pat \tilde \rho\|_3 \|\bfv\|_2\|\bfv\|_6  + \int_0^t \|\pat \tilde \rho\|_3\|\pat \bfomega_{\bfu}\times x\|_{\frac 83}^2\right)\\
\leq c\left(\int_0^t \|\pat \tilde \rho\|_3 \|\bfV\|_2\|\nabla \bfV\|_2 + c|\mbox{$\frac d{dt}$}\bfomega_{\mbox{\footnotesize{$\bfu$}}}|^2 \int_0^t \|\pat \tilde \rho\|_3 \right)\\
\leq  c(\varepsilon)\int_0^t\|\pat \tilde \rho\|_3^2\|\bfV\|_2^2 + \varepsilon \|\nabla \bfV\|_{L^2(L^2)}^2 + c\int_0^t\|\pat \tilde \rho\|_3 \int_\B |\mbox{$\frac d{dt}$}\bfomega_{\mbox{\footnotesize{$\bfu$}}}\times \bfx|^2\tilde\rho \\
\leq c(\varepsilon)\int_0^t\left(\|\pat \tilde \rho\|_3^2 + \|\pat \tilde \rho\|_3\right)\int_\S \tilde\rho |\bfW|^2 + \varepsilon \|\nabla \bfV\|_{L^2(L^2)}^2,
\end{multline}
 By a similar argument
\begin{multline}\label{odhad.prava.strana}
\int_0^t \int_\S \pat \bfF_3\cdot \bfW \leq \int_0^t \|\pat \bfF_3\|_{\frac 65} \|\bfV\|_{6} + c\int_0^t |\mbox{$\frac d{dt}$}\bfomega_{\mbox{\footnotesize{$\bfu$}}}|\|\pat \bfF_3\|_{\frac 65}\\
\leq c\int_0^t\|\pat \bfF_3\|^2_{\frac 65}  + \varepsilon \|\nabla \bfV\|^2_2 + \int_0^t \int_\S \tilde \rho |\bfW|^2
\end{multline}
and
\begin{equation}\label{odhad.prava.strana.1}
\int_0^t \int_\S \pat \bfF_4 :\nabla \bfW \leq c(\varepsilon)\int_0^t\int_\S |\pat \bfF_4|^2 + \varepsilon \int_0^t \int_\S |\nabla \bfV|^2 + c\int_0^t\int_\S \tilde \rho |\bfW|^2.
\end{equation}
From Korn's inequality we also get
$$
\int_0^t \int_\S S(\nabla \bfW):\nabla \bfW \geq c \int_0^t \int_\S |\nabla \bfV|^2.
$$
Thus, from (\ref{derivace.U}), (\ref{odhad.nabla.U}), (\ref{odhad.prava.strana}) and (\ref{odhad.prava.strana.1}) we deduce
\begin{multline}\label{pred.gronwallem}
\frac 12 \int_\S \tilde \rho |\bfW|^2(t) + c \int_0^t \int_\S |\nabla \bfV|^2 \\
\leq \frac 12 \int_\S \tilde \rho |\bfW|^2 (0) + \int_0^T \|\pat \bfF_3\|^2_{\frac 65} + \int_0^T\|\pat \bfF_4\|_{2}^2 + c(1 + \|\pat \tilde \rho\|^2_{L^2(L^3)})\int_0^t \int_\S \tilde \rho |\bfW|^2
\end{multline}
From (\ref{bolm.strong}) we deduce
\begin{multline*}
\int_\S \tilde \rho |\bfW|^2(0) = -\int_\S S(\nabla \bfu(0))\nabla :\bfW(0) + \int_\S \bfF(0) \bfW(0) = \int_\S \tilde\rho \bfW(0)\left(\diver S(\nabla \bfu(0))\frac 1{\tilde\rho}  + \bfF(0)\frac 1{\tilde \rho}\right) \\= \int_\S \diver S(\nabla \bfu(0))\left(\diver S(\nabla \bfu(0))\frac 1{\tilde\rho}  + \bfF(0)\frac 1{\tilde \rho}\right) + \int_\S \bfF(0)\left(\diver S(\nabla \bfu(0))\frac 1{\tilde\rho}  + \bfF(0)\frac 1{\tilde \rho}\right)\\
\leq c\left(\|\nabla^2 \bfw(0)\|_2^2 + \|\bfF(0)\|_2^2\right).
\end{multline*}
By the Gronwall inequality we get
\begin{multline*}
\|\bfW\|_{L^\infty(L^2)}^2 \leq c \left(\|\pat \bfF_3\|^2_{L^2(L^{\frac 65 })} + \|\pat \bfF_4\|^2_{L^2(L^2)} + \|\bfF(0)\|^2_{L^2}  + \| \bfw_0\|_{W^{2,2}}^2\right)\\\cdot \exp (c(1 + \|\pat \tilde \rho\|^2_{L^2(L^3)}))
\end{multline*}
and (\ref{pred.gronwallem}) also yields
\begin{multline*}
\|\bfV\|^2_{L^2(W^{1,2})}\leq c \left(\|\pat \bfF_3\|_{L^2(L^2)}^2  + \|\pat \bfF_4\|^2_{L^2(L^2)}+ \|\bfF(0)\|^2_{L^2} + \| \bfw_0\|_{W^{2,2}}^2\right)\\\cdot \left(1 + \|\pat \tilde \rho\|_{L^2(L^3)}^2\right) \exp(c(1 + \|\pat \tilde \rho\|_{L^2(L^3)}^2)).
\end{multline*}
In order to prove the remaining estimates, we again use classical regularity results for elliptic problems. Precisely, choosing $\bfphi \in W^{1,2}_0(\Ce)$ in (\ref{bolm.strong}) we get for a.a. $t$
$$
\int_\Ce S(\nabla \bfv):\nabla \bfphi = [\bfF,\bfphi] - \int_\Ce \tilde \rho \,\bfW \cdot\bfphi.
$$
which furnishes
$$
\|\nabla^2 \bfv\|^2_{L^\infty(L^2)} \leq c \left(\|\bfF_1\|_{L^\infty(L^2)}^2 + \|\bfF_2\|_{L^\infty(W^{1,2})}^2 + \|\tilde \rho \bfW\|_{L^\infty(L^2)}^2\right)
$$
and
$$
\|\nabla^3 \bfv\|^2_{L^2(L^2)} \leq c \left(\|\bfF_1\|_{L^2(W^{1,2})}^2  + \|\bfF_2\|_{L^2(W^{2,2})}^2+ \|\tilde \rho \bfW\|_{L^2(W^{1,2})}^2\right).
$$
Furthermore,
\begin{multline*}
\|\nabla \tilde \rho \bfW\|_{2}^2 = \int_\S |\nabla \tilde \rho|^2 |\bfW|^2 \leq \left(\int_\S |\nabla \tilde\rho|^6\right)^{\frac 13} \|\bfW\|_2 \|\bfW\|_6 \\
\leq \|\nabla\tilde \rho\|_6^2 \left(\|\bfW\|_2 \|\nabla\bfV\|_2 + \|\bfW\|_2 \|\mbox{$\frac d{dt}$}\bfomega_{\mbox{\footnotesize{$\bfu$}}} \times \bfx\|_6\right) \\ \leq \|\nabla \tilde \rho\|_6^4 \|\bfW\|_2^2 + \|\nabla \bfV\|_2^2 + \|\nabla \tilde \rho\|_6^2 \|\bfW\|_2^2
\end{multline*}
which, in turn, combined with (\ref{bolm.strong})$_3$ imply the last desired estimate.
\QED\end{proof}

\subsection{Proof of Theorem \ref{local.strong.solution}}
For $T>0$ we define a set
\begin{equation*}
\begin{split}
L_T = \{(\tilde \bfu, \tilde \sigma)&,\tilde \bfu \in C([0,T], W^{1,2}_R(\S)) \\
& \tilde \bfv \in C([0,T],W^{2,2}(\Ce)) \cap L^2(0,T,W^{3,2}(\Ce))\\
& \pat \tilde \bfu \in L^\infty(0,T,L^2_R(\S)) \cap L^2(0,T, W^{1,2}_R(\S))\\
&\tilde \sigma \in L^\infty(0,T,W^{2,2}(\Ce)),\ \pat \tilde \sigma \in L^\infty(0,T,W^{1,2}(\Ce))\\
&\|\tilde \bfu\|_{L^\infty(W^{1,2})}^2 + \|\tilde \bfv\|_{L^\infty(W^{2,2})}^2 + \|\tilde \bfv \|_{L^2(W^{3,2})}^2 + \|\pat \tilde \bfu \|_{L^\infty(L^2)}^2 + \|\pat \tilde \bfu\|_{L^2(W^{1,2})}^2\leq B_1\\&\|\tilde \sigma \|_{L^\infty(W^{2,2})}^2 + \|\pat \tilde \sigma\|_{L^\infty(W^{1,2})}^2 \leq B_2\\
& \tilde \bfu(0) = \bfu_0, \ \tilde \rho(0) = \rho_0 - \overline \rho, \ 0<\frac m2 \leq \tilde \sigma(t,x) +\overline \rho \leq 2M \},
\end{split}
\end{equation*}
where $B_1$ and $B_2$ are sufficiently large constants which will be chosen later.\\
We consider a mapping $\Phi$ defined on $L_T$ in the following way:
$$
\Phi(\tilde \bfu, \tilde \sigma) \mapsto (\bfu, \sigma),
$$
where $(\bfu, \sigma)$ is a solution to (\ref{bolm.strong}), (\ref{continuity.equation.strong}) with

\begin{equation*}
[\bfF,\bfphi] = \int_{\S}\big[(\tilde \rho\,\tilde \bfv \cdot \nabla\tilde \bfu + \tilde \rho \bfomega_{\mbox{\footnotesize{$\tilde\bfu$}}} \times \tilde\bfu)\cdot\bfphi + p(\tilde \rho) \diver\bfphi\big]\,, \ \
\tilde \rho = \tilde \sigma + \overline \rho\,\ \
\sigma_0 = \rho_0 - \overline \rho\,.
\end{equation*}
We also use the notation
$$
[\bfF,\bfphi] = \int_\S \underbrace{(\tilde \rho(\tilde \bfv \cdot \nabla)\tilde \bfu + \tilde \rho \bfomega_{\mbox{\footnotesize{$\tilde\bfu$}}} \times \tilde\bfu)}_{=\bfF_1}\cdot \bfphi + \int_\S \underbrace{p(\tilde \rho)}_{=\bfF_2} \diver \bfphi.
$$
and, by integration by parts,

$$
[\pat \bfF, \bfphi] = \int_\S  \underbrace{\left(\pat(\tilde \rho \bfomega_{\mbox{\footnotesize{$\tilde\bfu$}}} \times \tilde \bfu) + \pat (\tilde \rho\tilde \bfv) \nabla \bfu - \diver(\tilde \rho \tilde \bfv)\,\pat \tilde \bfu\right)}_{=\,{\displaystyle \pat} \bfF_3} \bfphi + \int_\S \underbrace{\left(\pat p(\tilde \rho){\bf1} {+} \tilde \rho \tilde \bfv \otimes\pat \tilde \bfu\right)}_{=\,{\displaystyle \pat} \bfF_4}: \nabla \bfphi.
$$
One may easily derive that
\begin{equation*}
\begin{split}
\|\bfF_1\|_{L^2(W^{1,2})}^2 &\leq T c(B_1,B_2)\,,\ \
\|\bfF_2\|_{L^2(W^{2,2})}^2  \leq T c(B_1,B_2)\\
\|\pat \bfF_3\|_{L^2(L^\frac 65)}^2 & \leq T c(B_1,B_2)\,,\ \
\|\pat \bfF_4\|_{L^2(L^2)}^2  \leq T c(B_1,B_2)\\
\|\nabla \tilde \rho\|_{L^4(L^6)}^4 + \|\pat \tilde \rho \|_{L^2(L^3)}^2 & \leq Tc(B_1,B_2)
\end{split}
\end{equation*}
Further, since  for every sufficiently smooth $f:(0,T)\mapsto \mathbb R_+$, it holds $\|f\|_\infty \leq |f(0)| + T^{\frac 12} \|f\|_2$, we get
\begin{equation*}
\begin{split}
\|\bfF_1\|_{L^\infty(L^2)}^2 & \leq T c(B_1,B_2) + \|\bfF_1(0)\|^2_{L^2}\\
\|\bfF_2\|_{L^\infty (W^{1,2})}^2 & \leq T c(B_1,B_2) + \|\bfF_2(0)\|^2_{W^{1,2}}\,,
\end{split}
\end{equation*}
while from Lemma \ref{continuity.equation.lemma} and Lemma \ref{bolm.strong.lemma} we infer
$$
\|\bfu\|_{L^\infty(W^{1,2})}^2 + \|\bfv\|_{L^\infty(W^{2,2})}^2 + \|\bfv \|_{L^2(W^{3,2})}^2 + \|\pat \bfu \|_{L^\infty(L^2)}^2\leq T c(B_1,B_2) + c\left(\|\bfF_1(0)\|_{L^2}^2 + \|\bfF_1(0)\|_{W^{1,2}}^2\right)
$$
and
$$
\|\sigma \|_{L^\infty(W^{2,2})} + \|\pat \sigma\|_{L^\infty(W^{1,2})} \leq (1+cB_1)(1+\|\sigma_0\|_{W^{2,2}})e^{T^{\frac 12} B_1}.
$$
We now choose $B_1$ sufficiently large, and $B_2$ such that
$$
B_2 > (1 + c B_1)(1 + \|\sigma_0\|_{W^{2,2}}).
$$
Furthermore,
$$
\|\sigma - \sigma_0\|_{L^\infty(W^{1,2})} \leq T\|\pat \sigma \|_{L^\infty(W^{1,2})} \leq Tc(B_1)
$$
and so, using well known interpolation results, we deduce
$$
\|\sigma - \sigma_0\|_{L^\infty((0,T)\times \Ce)} \leq c \|\sigma - \sigma_0\|_{L^\infty(W^{1,2})}^{\frac 13} \|\sigma - \sigma_0\|_{L^\infty(W^{2,2})}^{\frac 23} \leq c(B) T^{\frac 13} (1 + \|\sigma_0\|_2)^{\frac 23}.
$$
Consequently, there exists $T$ small enough such that
$$
\frac m2 \leq \sigma + \overline \rho \leq 2M \mbox{ a.e. in } \Ce\times [0,T].
$$
As a consequence, we may find some $T^*>0$ such that
$$
\Phi(L_{T^*}) \subset L_{T^*}.
$$
We next define the space $X := C([0,T^*], W^{1,2}_R(\S)) \times C([0,T^*], W^{1,2}(\Ce))$. It clearly follows that $L_{T^*}\subset X$ is closed and, from Arz\`ela-Ascoli theorem, also compactly embedded. Let us show that the function $\Phi$ is continuous on $L_{T^*}$. Indeed, suppose $\{(\tilde \bfu_n, \tilde \rho_n)\}_{n=1}^\infty\subset L_{T^*}$ is such that $(\tilde \bfu_n,\tilde \sigma_n)\rightarrow (\tilde \bfu, \tilde \sigma)$ in $X$. We subtract the equations for $(\bfu_n, \sigma_n)$ and $(\bfu, \sigma)$, multiply them by $(\bfu_n - \bfu, \sigma_n - \sigma)$ and, by the Gronwall lemma,  we get $( \bfu_n,\sigma_n)\rightarrow ( \bfu,\sigma)$ in $C([0,T], L^2(\S))$. Furthermore, using Gagliardo-Nirenberg inequality (see i.e. \cite[Theorem 1.24]{roub}) we get
$$
\|\nabla (\bfu_n - \bfu)\|_{2}\leq \|\bfu_n - \bfu\|^{\frac 56}_{W^{2,\frac {30}{23}}} \|\bfu_n - \bfu\|_2^{\frac 16}.
$$
Since $\bfu_n$ and $\bfu$ are from $L_{T^*}$ we get the uniform boundedness of $\|\bfu_n - \bfu\|^{\frac 56}_{W^{2,\frac {30}{23}}}$. The same argument applies
to the quantity $\sigma_n-\sigma$, which then allows us to deduce  $(\bfu_n, \sigma_n)\rightarrow (\bfu,\sigma)$ in $X$. As a result, 
all assumptions of Schauder theorem are fulfilled and we may thus conclude with the existence of  a fix point for the map $\Phi$, namely, a solution to (\ref{local.existence.strong}). 
As for the uniqueness of such a solution, its proof does not differ from the one presented in \cite[Section 3]{valli} and therefore will be omitted. Theorem \ref{local.strong.solution} is completely proved.
\setcounter{equation}{0}
\section{Global existence of strong solutions for small data}
Objective of this section is to show that the local in time solutions constructed in Theorem \ref{local.strong.solution}  can be extended to all times, provided the magnitude of the initial data is restricted in an appropriate sense. Before stating the main result, we introduce some {\em further notation}. Specifically, we denote by $[\ \cdot\ ]_k$ the  sum of $L^2$-norms involving only interior (in $\Ce$) and tangential derivatives (at $\partial\Ce$) of order $k$, and by $]| \cdot|[_k$, $[| \cdot|]_k$   suitable norms equivalent to the norm $\|\cdot\|_{k,2}$. Furthermore, we set
\begin{equation}\label{def.psi}
\psi(t) := ]|\bfv(t)|[_{1}^2 + ]|\sigma(t)|[_2^2 + [\bfv(t)]_2^2 + \int_\S \rho(t)|\pat \bfu(t)|^2 + \int_\S \rho(t) |\bfu(t)|^2 + \frac{p^\prime(\bar{\rho})}{\overline \rho} \|\pat \sigma(t)\|_2^2\,,
\end{equation}
where, we recall, $\sigma$ is defined in (\ref{s}). We also put
$$
\mathscr E(t):=\int_{\S}\rho\,|\bfu(t)|^2+\frac {2a}{\gamma-1}\int_{\Ce}\rho^\gamma(t)\,, 
$$ 
representing (twice) the total energy of the coupled system.
\smallskip\par
The following theorem holds.
\Bt\label{long.time.strong}
Let $ \Ce$ be of class $C^4$, $\bfu_0\in W^{1,2}_R(\S)$, $\bfu_0|_{\Ce}\in W^{2,2}(\Ce)$, $\rho_0\in W^{2,2}(\Ce)\cap L^\gamma(\Ce)$, $\gamma >1$. Then, there exists $\kappa_0>0$ such that if $\psi(0)+\mathscr E(0)\leq\kappa_0$, there are  uniquely determined $$\bfu\in C(\mathbb R^+, W^{1,2}_R(\S)),\ \mbox{with}\ \bfv \in C(\mathbb R^+, W^{2,2}(\Ce))\cap L^2_{\lokal}(\mathbb R^+, W^{3,2}(\Ce))$$ and
$$\rho\in C(\mathbb R^+, W^{2,2}(\Ce))$$
solving (\ref{first.sys.2}), (\ref{initcon}). Moreover,
$$
\pat \bfu\in C(\mathbb R^+, L^2_R(\S))\cap L^2_{\lokal}(\mathbb R^+, W^{1,2}_R(\S))\ \mbox{ and } \ \pat \rho \in C(\mathbb R^+, W^{1,2}(\Ce)).
$$
\Et
\Br
It is worth emphasizing that the ``smallness" condition required on the quantity $\psi(0)+\mathscr E(0)$ can be reformulated by imposing a restriction {\em only} on the size of the initial data $\bfu_0$ and $\rho_0$ in the norms of the function classes to which they belong. This statement is obvious for all terms appearing in the definition of $\psi$ and $\mathscr E$, with the exception of those involving the {\em time} derivative of $\bfu$ and $\sigma$. However, it is readily seen that the solution $(\bfu,\rho)$ of Theorem \ref{long.time.strong} possesses as much  regularity as to ensure the validity of  (\ref{first.sys.2})$_{1,3,4}$ (for a.a. $x\in\mathcal C$) and that of the derivative of (\ref{first.sys.2})$_{5}$ up to time $t=0$. Then, by employing standard Sobolev embedding theorems one can show the desired property. 
\Er
\par
The proof of Theorem \ref{long.time.strong} relies upon a number of global (in time) estimates that will be proved in the next sections.
\label{section.global.sol}
\subsection{Estimates of time derivatives}
We shall consider the following system
\begin{equation}\label{ge1}
\begin{split}
\int_\S \rho \,\pat \bfu\cdot \bfphi + \int_\Ce S(\nabla \bfu):\nabla \bfphi - \int_\Ce p_1 \rho \diver \bfphi &= [\bfh,\bfphi]\,,\,\mbox{ for all }\bfphi \in W^{1,2}_R(\S)\\
\pat \sigma + \bfv \nabla \sigma + \overline \rho \diver \bfv &= f_0,
\end{split}
\end{equation}
where $p_1 = \gamma\frac{p(\overline \rho)}{\overline \rho} (= p'(\overline \rho))$,
\begin{equation}\label{f0}
f_0 = -\sigma \diver \bfv,
\end{equation}
and
$$
[\bfh,\varphi] = \int_\S \bfh_1\cdot \bfphi + \int_\Ce \bfh_2\cdot \bfphi,
$$
with
\begin{equation}\label{h}
\begin{split}
\bfh_1 &= - \rho\, \omega_{\mbox{\footnotesize{$\bfu$}}}\times \bfu,\\
\bfh_2 &= -\rho \,\bfv\cdot \nabla\bfu + \left(p'(\overline\rho) - p'(\rho)\right) \nabla \sigma
\end{split}
\end{equation}
Furthermore, we assume that
\begin{equation}\label{bounded.density}
\frac{\overline \rho} 4 \leq \underbrace{\sigma(t,x) + \overline \rho}_{{\displaystyle=\rho}} \leq 3\overline \rho
\end{equation}
Since $\diver \bfphi = \diver(\bfphi - \bfomega_{\mbox{\footnotesize{$\bfphi$}}} \times \bfx - \bfxi_{\mbox{\footnotesize{$\bfphi$}}})$ and $\bfphi - \bfomega_{\mbox{\footnotesize{$\bfphi$}}} \times \bfx -\bfxi_{\mbox{\footnotesize{$\bfphi$}}}\in W^{1,2}_0(\Ce)$, we have
$$
\int_\Ce p_1 \rho \diver \bfphi = \int_\Ce p_1 \sigma \diver \bfphi.
$$
We choose $\bfphi = \bfu$ in (\ref{ge1})$_1$,  multiply (\ref{ge1})$_2$ by $\frac{p_1 \sigma}{\overline \rho}$ and add the resulting equations side by side. In view of (\ref{bounded.density}), we get
\begin{multline}\label{estimate.1}
\frac 12\frac d{dt} \int_\S \rho |\bfu|^2 + \frac{p_1}{\overline \rho}\frac d{dt} \|\sigma\|_2^2 + c\|\nabla \bfv\|_2^2 \leq \\ \varepsilon (\|v\|_{3,2}^2 + \|\sigma\|_2^2 + \|\bfu\|_2^2) + c(\varepsilon) (\|\bfh_1\|_2^2 + \|\bfh_2\|_{2}^2 +\|\sigma\|_2^4 + \|f_0\|_2^2 + \|\bfv \|_{2,2}^2\|\bfu\|_2^2),
\end{multline}
where we used
$$
\int_\Ce \diver \bfv |\sigma|^2 \leq c \|\diver \bfv\|_\infty \|\sigma\|_2^2 \leq c \|\bfv \|_{3,2} \|\sigma\|_2^2 \leq \varepsilon \|\bfv \|_{3,2}^2 + c(\varepsilon)\|\sigma\|_2^4
$$
and
\begin{multline*}
\left|\frac 12\int_\Ce \pat \rho |\bfu|^2\right| = \left|\int_\Ce \rho \bfv\cdot \nabla \bfu\cdot \bfu\right| \leq \int_\Ce \rho |\bfv| |\nabla \bfv| |\bfu| + c\int_\Ce \rho |\bfv| (|\bfomega_{\mbox{\footnotesize{$\bfu$}}}| + |\bfxi_{\mbox{\footnotesize{$\bfu$}}}|) |\bfu| \leq \\ c\left(\|\bfv\|_\infty \|\nabla \bfv\|_2\|\bfu\|_2 + \|\bfv\|_\infty (\|\bfomega_{\mbox{\footnotesize{$\bfu$}}}\|_2  + \|\bfxi_{\mbox{\footnotesize{$\bfu$}}}\|_2)\|\bfu\|_2\right)\leq c(\varepsilon) \|\bfv\|_{2,2}^2\|\bfu\|_2^2  + \varepsilon (\|\bfu\|_2^2 + \|\nabla \bfv\|_2^2).
\end{multline*}
Further, we differentiate (\ref{ge1}) with respect to $t$, (dot-)multiply the first equation by $\bfphi = \pat \bfu$ and the second equation by $\frac{p_1}{\overline \rho}\pat \sigma $. We get the following two equations
\begin{equation}
\begin{split}\label{time.estimate}
\frac 12 \int_\S \rho \pat |\pat \bfu|^2 + \int_\Ce S(\nabla \pat \bfv):\nabla \pat \bfv - \int_\Ce p_1 \pat \rho \diver \pat \bfv  + \int_\Ce \pat \sigma |\pat \bfu|^2&= [\pat \bfh,\pat\bfu]
\\
\frac 12 \frac{p_1}{\overline \rho} \int_\Ce \pat |\pat \sigma|^2 + \frac{p_1}{\overline \rho} \left(\int_\Ce  \pat \sigma\, \pat \bfv\cdot \nabla \sigma + \int_\Ce \pat \sigma\,\bfv\cdot \nabla \pat \sigma \right) + \int_\Ce p_1 \pat \sigma \diver \pat \bfv &= \int_\Ce \pat f_0 \frac{p_1}{\overline \rho} \pat \sigma.
\end{split}
\end{equation}
Recalling that $\pat \sigma = -\diver (\rho \bfv)$, we have
\begin{multline*}
\left|\int_\Ce \pat \sigma |\pat \bfu|^2\right| = \left|\int_\Ce \diver (\rho \bfv) |\pat \bfu|^2\right| = 2\left|\int_\Ce \rho \bfv\cdot \nabla \pat \bfu\cdot \pat \bfu\right|\\
\leq c \int_\Ce |\bfv| |\nabla \pat \bfu||\pat \bfu|\leq c \|\bfv\|_\infty \|\nabla \pat \bfu\|_2 \|\pat \bfu\|_2 \leq c\|\bfv\|_{2,2}^2 \|\pat \bfu\|_2^2 + \varepsilon \|\nabla \pat \bfu\|_2^2,
\end{multline*}
and
\begin{equation*}
\|\nabla \pat \bfu\|_2^2 = \|\nabla \pat \bfv + \nabla (\pat \omega \times x)\|_2^2 \leq c\left(\|\nabla \pat \bfv\|_2^2 + |\pat \omega|_2^2\right) \leq c\left(\|\nabla \pat \bfv\|_2^2 + \|\pat \bfu\|_2^2\right).
\end{equation*}
Consequently,
\begin{equation*}
\left|\int_\Ce \pat \sigma |\pat \bfu|^2\right|\leq c(\varepsilon)\|\bfv\|_{2,2}^2 \|\pat \bfu\|_2^2 + \varepsilon (\|\nabla \pat \bfv\|_2^2 + \|\pat \bfu\|_2^2).
\end{equation*}
We add together both relations in (\ref{time.estimate}) and use the same procedure as  in \cite[proof of Lemma 4.2]{valli} to get
\begin{multline}\label{estimate.2}
\frac 12\pat \int_\S\rho\pat |\bfu|^2 + \pat \frac{p_1}{\overline \rho}\|\pat \sigma\|_2^2 + c\|\nabla \pat \bfv\|_2^2 \leq \\
c(\varepsilon) \left(\left\|\pat \bfh_1\right\|_2^2 + \left\|\pat \bfh_2\right\|_{W^{-1,2}_R(\S)}^2 + \left\|\pat f_0\right\|_2^2 + \|\sigma\|_2^4 + \left\|\pat \sigma\right\|_2^4 + \|\bfv\|_{2,2}^2\|\pat \bfu\|_2^2\right) \\
+ \varepsilon \left(\|\pat \bfu \|_2^2 + \|\pat \sigma\|_2^2 + \|\bfv\|_{3,2}^2 \right)\,. 
\end{multline}
It remains to derive estimates for $\|\pat \bfu\|_2^2$. To this end, we take $\bfphi = \pat \bfu$ in (\ref{ge1})$_1$. We find
$$
\int_\S \rho |\pat \bfu|^2 + \int_\Ce S(\nabla \bfv):\nabla \pat \bfv - \int_\Ce p_1 \sigma \diver \pat \bfv = \int_\Ce \bfh_1\cdot \pat \bfu + \int_\C \bfh_2\cdot \pat \bfu
$$
which yields
$$
c\|\pat \bfu\|_2^2 \leq c \|\nabla \bfv\|_2\|\nabla \pat \bfv\|_2 + c\|\sigma\|_2\|\pat \nabla \bfv\|_2 + (\|\bfh_1\|_2+\|\bfh_2\|_2)\|\pat \bfu\|_2 .
$$
The latter in conjunction  with (\ref{bounded.density}), (\ref{estimate.1}) and (\ref{estimate.2}) furnishes
\begin{multline}\label{estimate.3}
\|\pat \bfu\|_2^2 \leq c \left(\|\bfh_1\|_2^2 + \|\bfh_2\|_2^2 + \|\sigma\|_{2}^4 + \|f_0\|_2^2 +\|\pat \bfh_1\|_2^2 + \|\pat \bfh_2\|_{W^{-1,2}_R(S)}^2 + \|\pat f_0\|_2^2 + \|\pat \sigma\|_2^4 \right.\\\left.+ \|\bfv \|_{2,2}^2\|\pat \bfu\|_2^2  + \|\bfv\|_{2,2}^2\|\bfu\|_2^2\right) + \varepsilon(\|\bfv\|_{3,2}^2 + \|\sigma\|_2^2 + \|\bfu\|_2^2 + \|\pat \sigma\|_2^2).
\end{multline}

\subsection{Further estimates of $\bfv$}
In this subsection we concern the following system (system for $\bfv$ divided by $\rho$)
\begin{equation*}
\begin{split}
\pat \bfv  + \frac 1{\overline \rho} \diver S(\nabla \bfv) + \sigma\,\nabla p_2   &= f\\
\bfv|_{\partial \Ce} &= 0\\
\bfv(0) &= \bfv_0
\end{split}
\end{equation*}
complemented with
\begin{equation*}
\pat \sigma +\bfv\cdot \nabla \sigma + \overline \rho \diver \bfv = f_0,
\end{equation*}
where $p_2 = \frac{p'(\overline\rho)}{\overline \rho}$, $f_0$ is defined in (\ref{f0}), and
$$
f = \frac{\sigma}{\overline \rho \rho} \diver S(\nabla \bfv) + \left(\frac{p'(\overline \rho)}{\overline \rho} - \frac{p'(\rho)}{\rho}\right)\nabla \sigma + \pat (\bfomega_{\mbox{\footnotesize{$\bfu$}}} \times \bfx) + \bfomega_{\mbox{\footnotesize{$\bfu$}}} \times \bfu + \bfv\cdot \nabla\bfu.
$$
Arguing exactly as in \cite[Section 4]{valli} we get the following estimates
\begin{equation}\label{valli.estimates}
\begin{split}
\frac{2\mu}{\overline\rho}\frac d{dt} \|\D (\bfv)\|_2^2 + \frac{\lambda - \frac 23 \mu}{\overline \rho}\frac d{dt} \|\diver \bfv\|_2^2  + c\|\bfv\|_{2,2}^2& \leq c\left(\|\sigma\|_{1,2}^2 + \|f\|_2^2\right)\\
\frac d{dt} \|\nabla \sigma\|_{1,2}^2&\leq c(\varepsilon) \left(\|\diver \bfv\|_{2,2}^2 + \|\sigma\|_{2,2}^2 + \|\sigma\|_{2,2}^4\right) + \varepsilon \|\bfv\|_{3,2}^2\\
\|\pat \sigma\|_{1,2}^2 &\leq c\left(\|\bfv \|_{2,2}^2 + \|f_0\|_{1,2}^2 + \|\bfv\|_{2,2}^2\|\sigma\|_{2,2}^2\right)\\
\|\sigma\|_{1,2}^2 + \|\bfv\|_{2,2}^2 &\leq c\left(\|\diver \bfv\|_{1,2}^2 + \|\pat \bfv\|_2^2 + \|f\|_2^2\right)\\
\|\sigma\|_{2,2}^2 + \|\bfv\|_{3,2}^2 &\leq c\left(\|\diver \bfv\|_{2,2}^2 + \|\nabla \pat \bfv\|_{2}^2 + \|f\|_{1,2}^2\right),
\end{split}
\end{equation}
\begin{multline}\label{valli.estimates2}
\frac d{dt} [\bfv]_1^2 + \pat |[\sigma]|_1^2 + \|\diver \bfv\|_{1,2}^2\\
\leq c(\delta)\left(\|\bfv\|_{1,2}^2 + \|\pat \bfv\|_2^2 + \|f_0\|_2^2 + \|f_0\|_{1,2}^2 + \|\sigma\|_{1,2}^4\right) + \delta \|\sigma\|_{1,2}^2 + \delta \|\bfv\|_{3,2}^2,
\end{multline}
and
\begin{multline}\label{valli.estimates3}
\frac d{dt} [\bfv]_{2}^2 + \frac d{dt} |[\sigma]|_2^2 + \|\diver \bfv\|_{2,2}^2\leq \\c(\varepsilon)\left(\|\bfv\|_{2,2}^2 + \|\pat \bfv\|_{1,2}^2 + \|\sigma\|_{1,2}^2 + \|f\|_{1,2}^2 + \|f_0\|_{2,2}^2 + \|\sigma\|_{2,2}^4\right) + \varepsilon \left(\|\sigma\|_{2,2}^2 + \|\bfv\|_{3,2}^2\right)
\end{multline}
where  $|[\ \cdot\ ]|_k$ is a norm equivalent to the $W^{k,2}(\Ce)$-norm. 

\subsection{Estimates of $\psi$}
We recall that the function $\psi$ is  defined in (\ref{def.psi}).
From (\ref{estimate.1}), (\ref{estimate.2}), (\ref{estimate.3}) and (\ref{valli.estimates}) we deduce
\begin{equation}\label{complicated1}
\begin{array}{l}\frac{2 \mu}{\overline\rho}\frac d{dt} \|\D(\bfv)\|_2^2 + \frac{\lambda - \frac 23 \mu}{\overline \rho}\pat \|\diver \bfv\|_2^2 + \pat \|\nabla \sigma\|_{1,2}^2 + \frac12\frac d{dt} \int_\B\rho|\pat \bfu|^2 + \frac{p_1}{\overline \rho}\frac d {dt} \|\pat \sigma\|_2^2 + \frac d{dt} \int_\B\rho|\bfu|^2 + \\ \\ \, \, +  \|\bfv\|_{3,2}^2 + \|\sigma\|_{2,2}^2 + \|\pat \bfv\|_{1,2}^2 + \|\pat \bfu\|_2^2 + \|\pat \sigma\|_{1,2}^2 \\\\\leq c(\varepsilon)\left(\|\sigma\|_{2,2}^2 + \|f\|_{1,2}^2 + \|\diver \bfv\|_{2,2}^2   + \|\sigma\|_{2,2}^4 + \|\pat \bfh_1\|_{2}^2 + \|\pat \bfh_2\|_{W^{-1,2}_R(\S)}^2 + \|\pat f_0\|_2^2 \right.\\\\ \left.+ \|\pat \sigma\|_{2}^4 + \|\bfv\|_{2,2}^2\|\pat \bfu\|_2^2 + \|\bfv\|_{2,2}^2 + \|f_0\|_{1,2}^2 + \|\bfv\|_{2,2}^2\|\sigma\|_{2,2}^2 + \|\bfv\|_{3,2}^2\|\bfu\|_2^2 + \|\bfh_1\|_2^2 + \|\bfh_2\|_2^2\right)\\\\ + \varepsilon(\|\sigma \|_{2}^2+\|\bfv\|_{3,2}^2 + \|\pat \bfu\|_2^2 + \|\pat \sigma\|_2^2 + \|\bfu\|_2^2).
\end{array} 
\end{equation}
We add to both sides of this inequality terms $\frac d{dt}[\bfv]_2^2$, $\frac d{dt} |[\sigma]|_2^2$, $\frac d{dt}[\bfv]_1^2$ and $\frac d{dt} |[\sigma]|_1^2$ multiplied by proper constants in order to use (\ref{valli.estimates2}) and (\ref{valli.estimates3}). Thus, (\ref{complicated1}) yields
\begin{multline}\label{complicated3}
\frac d {dt}\psi(t)
+ \|\bfv\|_{3,2}^2 + \|\sigma\|_{2,2}^2 + \|\pat \bfv\|_{1,2}^2 + \|\pat \bfu\|_2^2 + \|\pat \sigma\|_{1,2}^2 \\\leq c(\varepsilon)\left(\|f\|_{1,2}^2 +  \|\sigma\|_{2,2}^4 + \|\pat \bfh_1\|_{2}^2 + \|\pat \bfh_2\|_{W^{-1,2}_R(\S)}^2 + \|\pat f_0\|_2^2 + \|\pat \bfu\|_2^2\right.\\\left.+ \|\pat \sigma\|_{2}^4 + \|\bfv\|_{2,2}^2\|\pat \bfu\|_2^2  + \|f_0\|_{2,2}^2 + \|\bfv\|_{2,2}^2\|\sigma\|_{2,2}^2 + \|\bfv\|_{3,2}^2\|\bfu\|_2^2 + \|\bfh_1\|_2^2 + \|\bfh_2\|_2^2\right)\\ + \varepsilon(\|\bfv\|_{3,2}^2+\|\pat \bfu\|_2^2 + \|\pat \sigma\|_2^2 + \|\bfu\|_2^2 + \|\sigma\|_{2,2}^2)\,.
\end{multline}
Furthermore, taking also into account (\ref{h}), we show
\begin{equation}
\begin{split}\label{rhs.estimates}
\|f_0\|_{2,2}^2 \leq &c\|\sigma\|_{2,2}^2 \|\bfv\|_{3,2}^2\\
\|\pat f_0\|_2^2 \leq &c\left(\|\pat \sigma\|_2^2 \|\bfv\|_{3,2}^2 + \|\sigma\|_{2,2}^2 \|\pat \bfv\|_{1,2}^2\right)\\
\|f\|_{1,2}^2\leq &c\left(\|\sigma\|_{2,2}^2\|\bfv\|_{3,2}^2 + \|\sigma\|_{2,2}^4 + \|\bfu\|_2^4 + \|\bfu\|_2^2\|\bfv\|_{3,2}^2 + \|\bfv\|_{3,2}^2\|\bfv\|_{1,2}^2 + \|\pat \bfu\|_2^2\right)\\
\|\bfh_1\|_2^2 \leq &c\|\bfu\|_2^4\\
\|\pat \bfh_1\|_2^2\leq &c\left(\|\bfu\|_2^2\|\pat \bfu\|_2^2 + \|\pat\sigma\|_2^2\|\bfu\|_2^2\|\bfv\|_{3,2}^2\right)\\
\|\bfh_2\|_2^2\leq & c\left(\|\bfv\|_{3,2}^2\|\bfv\|_{1,2}^2 + \|\bfv\|_{3,2}^2\|\bfu\|_2^2 + \|\sigma\|_{2,2}^4\right).
\end{split}
\end{equation}
Observing that, by Sobolev inequality, $\|\cdot\|_{W^{-1,2}_R(\S)}\leq c\|\cdot\|_{\frac 65} \leq c\|\cdot\|_2$, we obtain
\begin{multline}\label{rhs.estimates2}
\|\pat \bfh_2\|_{W^{-1,2}_R(\S)}^2 \leq  c \left(\|\pat \bfv\|_2^2\|\bfv\|_{3,2}^2 + \|\pat \bfv\|_2^2\|\bfu\|_2^2 + \|\pat \sigma\|_2^2 \|\bfv\|_{1,2}^2\|\bfv\|_{3,2}^2 + \|\sigma\|_{2,2}^2\|\pat \sigma\|_{1,2}^2 \right.\\ \left.+ \|\pat \sigma \|_2^2 \|\bfu\|_2^2 \|\bfv\|_{3,2}^2\right).
\end{multline}
Combining (\ref{estimate.3}), (\ref{rhs.estimates})$_3$ and (\ref{rhs.estimates2}) we derive
\begin{multline}\label{rhs.estimate.2}
\|f\|_{1,2}^2\leq c\left(\|\sigma\|_{2,2}^2\|\bfv\|_{3,2}^2 + \|\sigma\|_{2,2}^4 + \|\bfu\|_2^4 + \|\bfu\|_2^2\|\bfv\|_{3,2}^2 + \|\bfv\|_{3,2}^2\|\bfv\|_{1,2}^2 + \|\bfu\|_2^2\|\pat \bfu\|_2^2\right.\\ + \|\bfv\|_{3,2}^2\|\pat\sigma\|_{2}^2 + \| \sigma\|_{2,2}^2\|\pat \bfv\|_{1,2}^2 + \|\sigma\|_{2,2}^2\|\pat\sigma\|_{1,2}^2\\\left. + \|\pat\sigma\|_2^2\|\bfv\|_{1,2}^2\|\bfv\|_{3,2}^2 + \|\pat\sigma\|_2^2\|\bfv\|_{3,2}^2\|\bfu\|_2^2+\|\pat\sigma\|_2^4 + \|\pat\bfu\|_2^2 \|\bfv\|_{3,2}^2 \right)\\
 + \varepsilon(\|\bfv\|_{3,2}^2 + \|\sigma\|_2^2 + \|\bfu\|_2^2 + \|\pat \sigma\|_2^2),
\end{multline}
where we used the inequality (see (\ref{10}))
$$
\|\pat \bfv\|_2^2 \leq c\left(\|\pat \bfu\|_2^2 +  |\frac d{dt}\bfomega_{\mbox{\footnotesize{$\bfu$}}}|^2\right)\leq c\|\pat \bfu\|_2^2.
$$
From (\ref{estimate.3}), (\ref{complicated3}), (\ref{rhs.estimates}), (\ref{rhs.estimates2}) and (\ref{rhs.estimate.2}), after absorbing some terms to left hand side --by taking $\varepsilon$ sufficiently small-- and adding $|\bfomega_{\mbox{\footnotesize{$\bfu$}}}|^2 + |\bfxi_{\mbox{\footnotesize{$\bfu$}}}|^2$ to both sides, we deduce 
\begin{equation}
\begin{split}\label{almost.final}
\frac d {dt} \psi + \|\bfv\|_{3,2}^2 &+ \|\sigma\|_{2,2}^2 + \|\pat \bfv\|_{1,2}^2 + \|\pat \bfu\|_2^2 + \|\pat \sigma\|_{1,2}^2  + |\bfomega_{\mbox{\footnotesize{$\bfu$}}}|^2 + |\bfxi_{\mbox{\footnotesize{$\bfu$}}}|^2\\
\leq &c\|\bfv\|_{3,2}^2\left(\|\pat \bfu\|_2^2 + \|\sigma\|_{2,2}^2 + \|\pat \sigma\|_2^2 + \|\bfv\|_{1,2}^2 + \|\bfu\|_2^2 + \|\pat \sigma \|_2^2\|\bfu\|_2^2 + \|\pat\sigma\|_2^2\|\bfv\|_{1,2}^2\right)\\
 + &c\|\sigma\|_{2,2}^2 \|\sigma\|_{2,2}^2 + c\|\pat \bfv\|_{1,2}^2\|\sigma\|_{2,2}^2
 + c\|\pat \sigma\|_{1,2}^2\left(\|\pat \sigma\|_{2}^2 + \|\sigma\|_{2,2}^2\right)
 + c\|\pat \bfu\|_2^2\|\bfu\|_2^2 + \\
 + &c\|\bfu\|_2^4 + |\bfomega_{\mbox{\footnotesize{$\bfu$}}}|^2 + |\bfxi_{\mbox{\footnotesize{$\bfu$}}}|^2.
\end{split}
\end{equation}
Since $\|\bfu\|_2\leq c(\|\bfv\|_{3,2} + |\bfomega_{\mbox{\footnotesize{$\bfu$}}}| + |\bfxi_{\mbox{\footnotesize{$\bfu$}}}|$),  from (\ref{almost.final}) we infer that there exist constants $\novak\label{c1}\geq 1$ and $\novak\label{c2}>0$ such that
\begin{equation}\label{derivace.psi}
\frac {d\psi}{dt} \leq -(\|\bfv\|_{3,2}^2 + \|\sigma\|_{2,2}^2 + \|\pat \bfv\|_{1,2}^2 + \|\pat \sigma\|_{1,2}^2  + \|\pat \bfu\|_2^2 + |\omega_\bfu|^2)(1-\starak{c1}(\psi + \psi^2)) + |\bfomega_{\mbox{\footnotesize{$\bfu$}}}|^2 + |\bfxi_{\mbox{\footnotesize{$\bfu$}}}|^2
\end{equation}
and
\begin{equation}\label{derivace.psi1}
\|\bfv\|_{3,2}^2 + \|\sigma\|_{2,2}^2 + \|\pat \bfv\|_{1,2}^2 + \|\pat \sigma\|_{1,2}^2 + \|\bfu\|_2^2 + \|\pat \bfu\|_2^2 + |\bfomega_{\mbox{\footnotesize{$\bfu$}}}|^2 + |\bfxi_{\mbox{\footnotesize{$\bfu$}}}|^2\geq \starak{c2} \psi.
\end{equation}
We next observe that, due to Sobolev's embedding theorem,  if $\psi(t)\leq \novak\label{c3}$ it follows that
\begin{equation}\label{sigma.estimate}
\frac{\overline \rho}2 \leq \sigma(t,x) + \overline \rho \leq \frac 32 \overline \rho.
\end{equation}
We then can show that there exists a constant $c$ such that
\begin{equation}\label{init.value.estimate}
\|\bfu(t)\|^2_{2} + \|\bfv(t)\|^2_{2,2} + \|\sigma(t)\|^2_{2,2} \leq c(\psi(t) + \psi^5(t)).
\end{equation}
Indeed, directly from (\ref{def.psi}) we have
$$
\|\bfu(t)\|^2_{2} + \|\sigma(t)\|_{2,2}^2 \leq c\psi(t).
$$
Moreover, from (\ref{ge1})$_1$ and regularity results for elliptic equation we have
$$
\|\bfv\|^2_{2,2} \leq \|\rho \pat \bfu\|^2_2 + \|\nabla p(\rho)\|^2_2 + \|\rho (\bfv \nabla) \bfu\|^2_2 + \|\bfomega_{\mbox{\footnotesize{$\bfu$}}}\times \bfu\|^2_2.
$$
Also,
$$
\|\bfomega_{\mbox{\footnotesize{$\bfu$}}}\times \bfu\|^2_2\leq c|\bfomega|^2\|\bfu\|^2_2 \leq c\|\bfu\|_2^4\leq c \psi^2 \leq c (1 + \psi^3)
$$
and, due to Sobolev embedding theorem and Gagliardo-Nirenberg interpolation theorem, 
\begin{multline*}
\|\bfv\cdot\nabla \bfv\|^2_2 \leq \|\bfv\|^2_\infty \|\bfv\|_{1,2}^2\leq \|\bfv\|_{1,4}^2\|v\|_{1,2}^2\leq c \|\bfv\|_{2,2}^{\frac 52} \|\bfv\|_6^{\frac 12} \|\bfv\|^2_{1,2} \leq c \|\bfv\|_{2,2}^{\frac 32} \|\bfv\|^{\frac 32}_{1,2}\\
\leq \varepsilon \|\bfv\|^2_{2,2} + c(\varepsilon) \|\bfv\|_{1,2}^6\leq \varepsilon \|\bfv\|^2_{2,2} + c(\varepsilon) \psi^5.
\end{multline*}
Therefore, after a suitable choice of $\varepsilon$, (\ref{init.value.estimate}) follows immediately.
\subsection{Proof of Theorem \ref{long.time.strong}}
We premise the following result.
\Bl
\label{boundedness.of.psi}
Let  $\psi=\psi(t)$ satisfy (\ref{derivace.psi}), (\ref{derivace.psi1}) for all $t\in [0,T)$. 
Then, there is $\kappa>0$ such that if 
\begin{equation}\label{SN}\psi(0) + \sup_{t\in[0,T]}\big(|\bfomega_{\mbox{\footnotesize{$\bfu$}}}(t)|+|\bfxi_{\mbox{\footnotesize{$\bfu$}}}(t)|\big)\leq \kappa\,,
\end{equation} necessarily $$\sup_{t\in [0,T]}\psi(t)\leq \kappa\,.$$
\El
\begin{proof} The proof is obtained by a slight modification of that given in  \cite[Lemma 4.10]{valli} and will be, therefore, omitted.
\QED\end{proof}
\par
In order to satisfy the hypothesis of this lemma, we only need to show an appropriate uniform bound on $|\bfomega_{\mbox{\footnotesize{$\bfu$}}}(t)|+|\bfxi_{\mbox{\footnotesize{$\bfu$}}}(t)|$. This will be obtained via the equation of energy balance for strong solutions that we may derive from (\ref{first.sys.1}). In fact, if we dot-multiply both sides of (\ref{first.sys.1})$_1$ by $\bfu$, integrate by parts over $\Ce\times (0,t)$, and take into account (\ref{first.sys.1})$_{2,3}$ and (\ref{pressure.rule}) (with $r\equiv\rho$), we deduce
\begin{multline*}
\frac12\int_{\Ce}\big[ \rho(t)|\bfu(t)|^2 + \frac{2p(\rho(t))}{\gamma -1}\big] + \int_{0}^t \int_{\mathcal C} S(\nabla\bfv(s)):\nabla \bfv(s) =  \frac12\int_{\mathcal C}\big[ \rho_0|\bfu_0|^2 + \frac{2p(\rho_0)}{\gamma -1}\big]\\ + \int_0^t\Big[\bfomega\cdot\int_{\partial\Ce}\bfx\times T(\bfu(s),p(\rho(s))\cdot\bfn +\bfxi\cdot\int_{\partial\Ce} T(\bfu(s),p(\rho(s))\cdot\bfn\Big]\,.
\end{multline*}
If we now employ (\ref{first.sys.1})$_4$ and (\ref{RM}) in the last two integrals of these equality, we conclude, for all $t\geq0$,
\begin{equation}\label{weneq}\begin{split}
\mathscr E(t)+2\int_0^t\int_{\Ce}S(\nabla\bfv(s))&:\nabla\bfv(s)=\mathscr E(0)\,,\\ \mathscr E:=\|\sqrt{\rho}\bfu\|_2^2+\frac2{\gamma-1}\|\sqrt{p(\rho)}\|_2^2&+\bfomega_{\mbox{\footnotesize{$\bfu$}}}\cdot\bfI_C\cdot\bfomega_{\mbox{\footnotesize{$\bfu$}}}+m_{\B}|\bfxi_{\mbox{\footnotesize{$\bfu$}}}|^2\equiv \int_{\S}\rho\,|\bfu|^2+\frac{2a}{\gamma-1}\int_{\Ce}\rho^\gamma\,.\end{split}
\end{equation}
This equation furnishes, in particular, the following uniform bound
\begin{equation}\label{DP}
\sup_{t\ge 0}\big(|\bfomega_{\mbox{\footnotesize{$\bfu$}}}(t)|+|\bfxi_{\mbox{\footnotesize{$\bfu$}}}(t)|\big)\leq c_4\,\mathscr E(0)\,,
\end{equation}
from which the validity of the hypothesis of the lemma follows by taking $\mathscr E(0)$ less than a suitable constant.

We are now in a position to prove Theorem \ref{long.time.strong}. According to Theorem \ref{local.strong.solution} there exist a time $T^*$ and functions $\bfu, \rho$ solving (\ref{local.existence.strong}) in $(0,T^*)$ in the function class there specified. We assume that $\psi(0)+\mathscr E(0)$  is sufficiently small, so that, by (\ref{DP}), condition (\ref{SN}) is satisfied. As a consequence, we may take $T^*$ depending only on the physical parameters and on $\kappa$ (see \cite[Section 4]{valli} for details). Then, by Lemma \ref{boundedness.of.psi}  and (\ref{weneq}), we get 
\begin{equation}\label{long.time.psi.bound}\psi(t)+\big(|\bfomega_{\mbox{\footnotesize{$\bfu$}}}(t)|+|\bfxi_{\mbox{\footnotesize{$\bfu$}}}(t)|\big)\leq \kappa\,, 
\end{equation}
for all $t\in[0,T^*)$. Furthermore, by (\ref{sigma.estimate}) and (\ref{init.value.estimate}) we have
$$\|\bfu(T^*)\|^2_2 + \|\bfv(T^*)\|^2_{2,2} + \|\sigma(t)\|^2_{2,2}\leq c_0
$$
for some constant $c_0>0$, and thus we may again use Theorem \ref{local.strong.solution} to establish existence of a (unique) solution on a time interval $(0,2T^*)$. We repeat this argument on each interval $(0,nT^*)$, $n\in \mathbb N$  in order to get the existence of solution on the whole positive real line. Theorem \ref{local.strong.solution}, (\ref{def.psi}) and (\ref{init.value.estimate}) yield that $\rho$ and $\bfu$ belong to the desired spaces.
\setcounter{equation}{0}
\section{Long time behavior of strong solutions}
Throughout this section, we shall work with strong solutions  constructed by Theorem \ref{long.time.strong}. We would like to point out that while, in general, weak solutions may admit multiple zero-velocity-limit solutions \cite{FePe},  the same issue does not happen in our case since the density of strong solutions is always bounded away from zero. 
\smallskip\par
We begin with the following simple observation that we state in the form of a lemma.

\Bl\label{infinity.limit}
Let $f\in C(\mathbb R^+)$, $f\geq 0$ be such that
$$
\int_0^\infty f(t){\rm d}t = c<\infty,\quad |f'(t)|\leq d<\infty\ \mbox{for all}\ t\in \mathbb R^+.
$$
Then $\lim_{t\to\infty} f(t)=0$.
\El
\begin{proof}
The claim follows at once from the following identity
$$
|f^2(t)-f^2(s)|=2\left|\int_s^tf(\tau)\,f^\prime(\tau)\right|\,,\ \ \mbox{all $t,s>0$}\,,
$$
and the assumptions.
\QED\end{proof}

We next give a definition of $\Omega$-limit set, appropriate to our problem.  
\Bd
Let $(\bfv, \bfomega_{\mbox{\footnotesize{$\bfu$}}}, \bfxi_{\mbox{\footnotesize{$\bfu$}}}, \rho)$ be a solution constructed in Theorem \ref{long.time.strong}. The corresponding $\Omega$-limit set, $\Omega(\bfv, \bfomega_{\mbox{\footnotesize{$\bfu$}}}, \bfxi_{\mbox{\footnotesize{$\bfu$}}}, \rho)\subset L^2(\Ce)\times \mathbb R^3\times \mathbb R^3 \times L^2(\Ce)$ is the set of all $(\hat \bfv, \hat \bfomega, \hat \bfxi, \hat \rho)$, for which there exists an increasing,  unbounded sequence $\{t_n\}\subset (0,\infty)$ such that
$$
\lim_{n\rightarrow \infty} (\|\bfv(t_n) - \hat \bfv\|_2 + |\bfomega_{\mbox{\footnotesize{$\bfu$}}}(t_n) - \hat \bfomega| + |\bfxi_{\mbox{\footnotesize{$\bfu$}}}(t_n) - \hat\bfxi| + \|\rho(t_n) - \hat \rho\|_2) = 0.
$$
\Ed
\Bl The above defined $\Omega-$limit set possesses the following properties:\label{prop.omega}
\begin{itemize}
\item[{\rm (i)}] It is not empty, compact and connected.
\item[{\rm (ii)}] Every element of the set is of the form $({\bf 0},\hat \bfomega, \hat \bfxi, \hat \rho)$ with $\hat \rho \in (\overline \rho/2, 3/2\overline\rho)$.
\item[{\rm (iii)}] It is invariant under solutions constructed in Theorem \ref{long.time.strong}.
\end{itemize}
\El
\begin{proof}
The first property follows from the uniform estimates proved  in Theorem \ref{long.time.strong} and basic properties of $\Omega-$limit sets (see e.g. \cite[Lemma 3]{Hal})
The second one is a consequence of Lemma \ref{infinity.limit}. Indeed, from the energy equality (\ref{weneq}) along with the Poincar\'e inequality we infer $\int_0^\infty\|\bfv\|^2_2\leq c<\infty$. Furthermore, (\ref{first.sys.1})$_1$ together with regularity of the constructed solution yields $|\pat\|\bfv(t)\|_2^2|\leq c$, uniformly in $t\geq0$, so that assumptions of Lemma \ref{infinity.limit} are verified and we conclude with the desired property. We point out that the smoothness of strong solution together with the convergence of $\bfv$ also yield
\begin{equation}\label{convergence.bv}
\|\bfv(t)\|_{1,2}\to 0\ \mbox{as}\ t\to \infty.
\end{equation}
The property $\overline\rho/2<\hat\rho<3/2\overline\rho$ is true due to (\ref{sigma.estimate}) and (\ref{long.time.psi.bound}), the latter valid for all $t\geq 0$. In order to show  the last property (iii), 
 let ${\sf s}(t):=(\bfv_\infty, \bfomega_\infty, \bfxi_\infty, \rho_\infty)$ denote the solution constructed in Theorem \ref{long.time.strong} emanating from $(0,\hat \bfomega,\hat \bfxi, \hat \rho)\in \Omega(\bfv, \bfomega_{\mbox{\footnotesize{$\bfu$}}}, \bfxi_{\mbox{\footnotesize{$\bfu$}}}, \rho)$.~\footnote{It's worth pointing out that, by definition of $\Omega$ and Theorem \ref{long.time.strong}, every point in $\Omega(\bfv, \bfomega_{\mbox{\footnotesize{$\bfu$}}}, \bfxi_{\mbox{\footnotesize{$\bfu$}}}, \rho)$ satisfies the assumptions on the data of that theorem, and thus can be used as initial condition.} We have to prove   that, for every $t>0$, ${\sf s}(t)\in \Omega(\bfv, \bfomega_{\mbox{\footnotesize{$\bfu$}}}, \bfxi_{\mbox{\footnotesize{$\bfu$}}}, \rho)$. Denote by $t_n$ the sequence satisfying $(\rho(t_n),\bfomega_{\mbox{\footnotesize{$\bfu$}}}(t_n), \bfxi_{\mbox{\footnotesize{$\bfu$}}}(t_n))\rightarrow (\hat\rho, \hat \bfomega, \hat \bfxi)$ in $L^2\times\mathbb R^3\times\mathbb R^3$. Since the strong solution is unique, it is enough to show that
$$
\bfv_n(t) \rightarrow \bfv_\infty(t) \mbox{ in } L^2, \quad
\rho_n(t)\rightarrow \rho_\infty(t) \mbox{ in }L^2, \quad \bfomega_n(t)\rightarrow \bfomega_\infty(t)\  \mbox{and}\ \bfxi_n(t)\rightarrow \bfxi_\infty(t)\mbox{ in }\mathbb R^3
$$
where $(\bfv_n(t),\rho_n(t),\bfomega_n(t), \bfxi_n(t))$ is a solution with initial conditions $(\bfv(t_n),\rho(t_n),\bfomega(t_n), \bfxi(t_n))$.%
Clearly, both $(\bfu_\infty, \rho_\infty)$  and $(\bfu_n, \rho_n)$ satisfy (\ref{first.sys.1}). If we now subtract, side by side,
We subtract (\ref{first.sys.1})$_1$ written for $(\bfu_\infty, \rho_\infty)$ from the one for $(\bfu_n,\ \rho_n)$. We next multiply the resulting equation by $\bfu_n -\bfu_\infty$. Furthermore, we  perform the same procedure for (\ref{first.sys.1})$_2$ and  multiply the resulting equation by $\rho_n - \rho_\infty$. After a  cumbersome but straightforward calculation that employs also the bounds proved in Theorem \ref{long.time.strong} we get for each $t>0$
\begin{multline*}
\frac 12\int_0^t \pas\left(\int_{\S}\rho_n |\bfu_n - \bfu_\infty|^2\right)  + \frac12\int_0^t \pas \|\rho_n -\rho_\infty\|_2^2 + \int_0^t \|\nabla(\bfv_n - \bfv_\infty)\|_2^2\leq\\
c\int_0^t\|\rho_n - \rho_\infty\|_2^2 + c\int_0^t\|\bfu_n - \bfu_\infty\|_2^2 + \varepsilon_0^t\int_0^t\|\nabla(\bfv_n - \bfv_\infty)\|_2^2  \\+ c\int_0^t\|\nabla \bfv_n\|_\infty^2 \|\rho_n - \rho_\infty\|_2^2 + c\int_0^t\|\pas \bfu_\infty\|_6^2\|\rho_n-\rho_\infty\|_2^2 + c\int_0^t\|\nabla \bfv_n\|_2^2 + c\int_0^t\|\nabla \bfv_n\|_2.
\end{multline*}
By Gronwall inequality
\begin{multline*}
\|(\bfu_n - \bfu_\infty)(t)\|_2^2 + \|(\rho_n - \rho_\infty)(t)\|_2^2\leq \\
c\big[\|(\bfu_n - \hat\bfu)(0)\|_2^2 + \|(\rho_n - \hat\rho)(0)\|_2^2 + \int_0^t\|\nabla \bfv_n\|_2^2 + \int_0^t\|\nabla \bfv_n\|_2\big]e^{\int_0^t\left( 1+\|\nabla \bfv_n\|_\infty^2 + \|\pas \bfu_\infty\|_6^2\right)}.
\end{multline*}
The desired property follows by letting $n\to\infty$ in this inequality and using (\ref{convergence.bv}), estimates of Theorem \ref{long.time.strong} and the fact that $(\bfu_n(0), \rho_n(0))\rightarrow (\hat{\bfu},\hat{\rho})$ strongly in $L^2$.
\QED\end{proof}
The next result shows that $\Omega(\bfv, \bfomega_{\mbox{\footnotesize{$\bfu$}}}, \bfxi_{\mbox{\footnotesize{$\bfu$}}}, \rho)$ is a subset of the set  of steady-state solutions.
\Bl\label{omega.stationary}
For any solution of Theorem \ref{long.time.strong}, the generic element of the corresponding set $\Omega(\bfv, \bfomega_{\mbox{\footnotesize{$\bfu$}}}, \bfxi_{\mbox{\footnotesize{$\bfu$}}}, \rho)$ is of the type $(\bfv\equiv{\bf 0},\rho_s,\bfomega_s,\bfxi_s)$ with $(\rho_s,\bfomega_s,\bfxi_s)$ satisfying  (\ref{stationary2}).  Moreover,
\begin{equation}\label{anmo}
|\bfI(\rho_s)\cdot\bfomega_s|=M_0\,,\ \ \int_{\Ce}\rho_s=m_{\mathcal F}\,,
\end{equation}
where $\bfI(\rho_s)$ is defined in (\ref{I}),  $M_0$ is the magnitude of the initial angular momentum $\bfM$ defined in (\ref{first.sys.2})$_4$, and, we recall, $m_{\mathcal F}$ is the mass of the fluid. 
\El
\begin{proof}
Due to the invariance of the $\Omega$-limit set, we derive from (\ref{first.sys.1}) that every solution $({\bf 0},\bfomega,\bfxi, \rho)$ emanating from a  point in $\Omega$ satisfies ($\ \dot{}\equiv \frac d{dt}$)
\begin{equation}\label{w.limit.sys}
\begin{split}
\partial_t\,{\rho} & = 0,\\
\rho(\dot \bfomega\times \bfx + \dot \bfxi) + \rho \bfomega\times(\bfomega\times \bfx + \bfxi) + \nabla  p(\rho) &= {\bf 0}.\\
m_\B\,\bfxi& =-\bfomega\times \int_{\Ce}\rho\bfx\,,
\end{split}
\end{equation}
and, moreover, $\rho\geq\nu>0$. The first of the above equations gives that $\rho$ is at most a function of $x$. We next apply the operator $\nabla \times $ on both sides of (\ref{w.limit.sys})$_{2}$ and obtain
$$
-2 \rho\,\dot \bfomega + \nabla \rho \times (\dot\bfomega \times \bfx+\dot\bfxi +\bfomega\times(\bfomega\times\bfx+\bfxi)) = {\bf 0}\ \mbox{ in }\Ce.
$$
Employing (\ref{w.limit.sys})$_2$ in the latter and recalling (\ref{pressure.rule}) (with $r\equiv\rho$), we deduce
\begin{equation}\label{stability.of.angular.momenta}
2 \dot \bfomega = -\frac 1{\rho} \nabla \rho \times (\frac1\rho\nabla p(\rho))=-a\,\frac\gamma\rho\nabla\rho\times\big(\rho^{\gamma-2}\nabla\rho\big)={\bf 0}.
\end{equation}
Thus (\ref{stability.of.angular.momenta}) yields $\dot\bfomega = {\bf 0}$, which, once combined with (\ref{w.limit.sys})$_3$ furnishes also $\dot\bfxi={\bf 0}$. We thus get that $\rho$, $\bfomega$ and $\bfxi$ must satisfy (\ref{sasi}) and, therefore, (\ref{stationary2}). It remains to show (\ref{anmo}). Condition (\ref{anmo})$_2$ follows by letting $t\to\infty$ (along a sequence) in (\ref{conmass}).  As for the first condition, we observe that, from (\ref{first.sys.2})$_4$ we get
\begin{equation}\label{M_0_0}
|\bfI_C\cdot\bfomega(t)+\int_{\Ce}\rho(t)\bfx\times\big(\bfv(t)+\bfomega(t)\times\bfx+\bfxi(t)\big)|=M_0\,,\ \ \mbox{for all $t\ge 0$}.
\end{equation} 
The claimed property is then shown by letting $t\to\infty$ (along a sequence) in this relation, and employing (\ref{convergence.bv}) and definition (\ref{I}).
\QED\end{proof}

\Br\label{bu} Since the $\Omega$-limit set is not empty, by the previous lemma we deduce that the set of solutions to (\ref{stationary2}) (or, equivalently, weak solutions to (\ref{stationary})) is not empty as well.
\Er
In view of the lemma just proved, the natural question is whether and when the generic solution of Theorem \ref{long.time.strong} will tend to a uniquely determined steady-state, namely, the corresponding $\Omega$-limit set reduces to a singleton. It is easy to show that this indeed happens in the (trivial and non-generic) case  $M_0=0$. In fact, we have the following result.

\Bt\label{M_0} Let $M_0=0$ and let $(\bfu,\rho)$ be a generic solution constructed in Theorem \ref{long.time.strong}. Then
$$
\Omega(\bfv, \bfomega_{\mbox{\footnotesize{$\bfu$}}}, \bfxi_{\mbox{\footnotesize{$\bfu$}}}, \rho)=\{{\bf 0}, {\bf 0}, {\bf 0}, m_{\F}/|{\Ce}|\}\,.
$$
\Et
\begin{proof} We observe that, in view of Lemma \ref{stab}, from (\ref{anmo})$_1$ it follows at once $\bfomega_s={\bf 0}$, which, once replaced in (\ref{sasi}), furnishes $\bfxi_s={\bf 0}$. Finally, from (\ref{stationary2}) and (\ref{anmo})$_2$ we conclude $\rho_s=\textrm{const.}=m_{\F}/|\Ce|$. 
\QED\end{proof}

\Br
Theorem \ref{M_0} is (for small initial data) the compressible counterpart of the same result shown in \cite{ST2} in the incompressible case.
\Er

If $M_0\neq 0$, the situation appears to be more complicated, and we are able to give an answer only under suitable assumption on the ``mass distribution" of the coupled system, and requiring the material constant $a$ to be ``sufficiently large". More precisely, in the following Lemma \ref{isolated.solutions} we show that if the tensor $\bfI(\bar{\rho})$, with $\bar{\rho}:=m_\F/|\Ce|$ and $\bfI(\rho)$ defined in (\ref{I}), has distinct eigenvalues, then all possible solutions to the system (\ref{u}), (\ref{stationary2}), (\ref{anmo}) must be isolated for $a$ large enough. (Notice that, by Lemma \ref{stab}, $\bfI(\bar{\rho})$ possesses three positive eigenvalues, for any $\bar{\rho}>0$.) As a consequence, since $\Omega$ is connected, it must coincide with one of these solutions; see Theorem \ref{final.claim}.
\Bl\label{isolated.solutions}
Suppose the eigenvalues of $\bfI(\bar\rho)$ are distinct. Then, there exists $a_0>0$ such that if $a>a_0$, the $\Omega-$limit set associated to a generic solution of Theorem  \ref{long.time.strong} reduces to a singleton.
\El
\begin{proof}
We recall that, by Lemma \ref{omega.stationary}, the generic element of the $\Omega-$limit set $(\rho_s,\bfomega_s,\bfxi_s)$ satisfies the following set of equations
\begin{equation}\label{adieu}
\begin{split}
\rho_s= \left(\frac{\gamma -1}{2a\gamma}(|\bfomega_s\times \bfx|^2 - (\bfomega_s\times \bfxi_s)\cdot\bfx) + c_s\right)^{\frac 1{\gamma - 1}}&=:{\sf f}(\omega_s,\xi_s, c_s),
\\
\bfomega_s\times (\bfI(\rho_s)\cdot \bfomega_s)& = 0,\\
|\bfI(\rho_s) \cdot \bfomega_s| &= M_0,\\
m_\S \bfxi_s& = \bfg(\rho_s)\times \bfomega_s\\
\int_{\Omega} {\sf f}(\bfomega_s, \bfxi_s, c_s)& = m_\F\,,
\end{split}
\end{equation}
where, in view of Theorem \ref{M_0}, we may assume  $M_0\neq 0$.
From the second equation above, we deduce that $\bfomega_s$ is an eigenvector of $\bfI(\rho_s)$ corresponding to some (positive) eigenvalue $\lambda_s$.  We  claim that if $\lambda_s$ is simple, then the solution $(\rho_s,\bfomega_s,\bfxi_s)$ is isolated. To this end, we notice that (\ref{adieu}) can be written in the following form
\begin{equation}\label{readytoimplicit}
\begin{split}
\left(\lambda_s{\bf 1} - \bfI({\sf f})\right)\cdot\bfomega_s & = 0,\\
\lambda_s^2 |\bfomega_s|^2 - M_0^2 & = 0,\\
m_\S \bfxi_s  - \bfg({\sf f})\times \bfomega_s & = 0,\\
\int_{\Ce} {\sf f}(\omega_s, \xi_s, c_s)& = m_\F,
\end{split}
\end{equation}
in the unknowns $\lambda_s\in \mathbb R^+$, $\bfomega_s, \bfxi_s \in \mathbb R^3$ and $c_s\in \mathbb R^+$. 
Consider the function $F:\mathbb R^8\mapsto \mathbb R^8$ 
defined as follows
$$
F(\lambda, \bfomega, \bfxi,c) = \left(\left(\lambda{\bf 1} - \bfI\right)\cdot\bfomega, \lambda^2 |\bfomega |^2 - M_0^2,m_\S \bfxi  
- 
\bfg
\times \bfomega, \int_\Ce {\sf f}(\bfomega,\bfxi, c)-m_\F\right)\,,
$$
where both $\bfI$ and $\bfg$ are evaluated at ${\sf f}(\bfomega,\bfxi,c)$ 
.
Our claim will then become a consequence of the implicit function theorem, provided we show
that $\nabla F(\lambda_s,\bfomega_s,\bfxi_s,c_s)$ is a regular matrix.
Let 
\begin{equation}\label{pi}
{\sf p}(\bfomega,\bfxi, c) := \frac{\gamma -1}{2a\gamma}\left(|\bfomega\times \bfx|^2 - 2(\bfomega\times \bfxi)\cdot\bfx\right) + c.
\end{equation}
From  (\ref{sigma.estimate}) it follows that there is $C=C(\overline \rho,\gamma)>0$ such that
\begin{equation}\label{chsf}
C>{\sf p}(\bfomega_s,\bfxi_s,c_s)>C^{-1}\,.
\end{equation}
Next, by a straightforward calculation we show 
$$
\nabla F =\left(
\begin{matrix}
\bfomega & \lambda{\bf1} - \bfI & 0& 0\\
2\lambda |\bfomega|^2 & 2\lambda^2 \bfomega & 0&0\\
0 & \mathbb S(\bfg)& m_\S {\bf 1}& 0\\
0 & 0 & 0 & \frac1{\gamma-1}\int_\Ce  {\sf p}^{\frac{2-\gamma}{\gamma - 1}}
\end{matrix}
\right) + \bfN(\bfomega,\bfxi,c)
$$
where $\mathbb S$ is defined in (\ref{S1}) and
$$\bfN(\bfomega,\bfxi,c) =\left(\begin{matrix}
0&-\frac{\partial\bfI}{\partial \bfomega}\cdot\bfomega & - \frac{\partial\bfI}{\partial \bfxi}\cdot\bfomega&0\\
0&0&0&0\\
0&-\frac{\partial \bfg}{\partial \bfomega}\times \bfomega&-\frac{\partial \bfg}{\partial \bfxi}\times \bfomega&0\\
0& \int_\Ce\frac{\partial {\sf f}(\bfomega,\bfxi,c)}{\partial \bfomega} & \int_\Ce\frac{\partial {\sf f}(\bfomega,\bfxi,c)}{\partial \bfxi} & 0
\end{matrix}\right). $$
Let us prove that there exists $C_0$ depending at most on $\bfxi_s,\ \bfomega_s,\ c_s,\ \mathcal C$ and $\gamma$ such that
\begin{equation}\label{boundonN}
|\bfN({\bfomega_s,\bfxi_s,c_s})|\leq a^{-1/(\gamma - 1)}C_0.
\end{equation}
We show the validity of this bound only for one entry of $\bfN$, as the proof for all other entries is exactly the same (or simpler). So, for example, we have
\begin{equation*}
\frac{\partial \bfI}{\partial \bfomega} = \frac{\partial\bfI}{\partial {\sf f}}\cdot\frac{\partial {\sf f}}{\partial \bfomega}= a^{-1/(\gamma-1)}\,{\sf p}^{\frac{\gamma-2}{\gamma-1}}(\bfomega,\bfxi,c)\,\bfK(\gamma,\mathcal C)\cdot\left[\left(2(\bfomega\times\bfx)\cdot\mathbb S(\bfx)-2\mathbb S(\bfxi)\cdot\bfx\right)\right]\,, 
\end{equation*}
with $\bfK(\gamma,\Ce)$  a tensor depending only on $\gamma$ and $\Ce$. Therefore,
 (\ref{boundonN}) follows, once we take into account (\ref{chsf}).
We next prove that the matrix $\nabla F - \bfN$ at the  point $(\lambda_s,\omega_s,\xi_s)$, denoted by $\bfM_s$, is regular for all values of $a>0$ and, moreover, also in the limit $a\to\infty$. Indeed, the only term depending on $a$ is the entry involving ${\sf p}$ which is (of course) defined for all $a>0$ and, in view of (\ref{adieu})$_{1,4}$ and (\ref{pi}) tends to a non-zero constant in the limit $a\to\infty$. Thus, to show that $\bfM_s$ is non-singular it is sufficient to prove that the reduced matrix
$$
{\bfP}:=\left(\begin{matrix}
 \bfomega_s & \lambda_s{\bf 1} - \bfI(\rho_s)\\ 2\lambda_s|\bfomega_s|^2 & 2\lambda_s^2\bfomega_s
\end{matrix}
\right)
$$
spans (uniquely) the whole of $\mathbb R^4$. But this is true since, denoting by $ \lambda_s,\ \lambda_2,\ \lambda_3$ the eigenvalues of $\bfI(\rho_s)$ and ${\bfr}_1=\bfomega_s,\ {\bfr}_2,\ {\bfr}_3$  corresponding eigenvectors, we get
\begin{equation*}
\begin{split}
{\bfP}\cdot\left(\begin{matrix} 0\\{\bfr}_2 \end{matrix}\right) & = (\lambda_s - \lambda_2)\left(\begin{matrix}{\bfr}_2\\0\end{matrix}\right),\\
{\bfP}\cdot\left(\begin{matrix} 0\\{\bfr}_3 \end{matrix}\right) & = (\lambda_s - \lambda_3)\left(\begin{matrix}{\bfr}_3\\0\end{matrix}\right),\\
{\bfP}\cdot\left(\begin{matrix} 0\\{\bfr}_1 \end{matrix}\right) & = \left(\begin{matrix}0\\2\lambda_s^2|{\bfr}_1|^2\end{matrix}\right),\\
{\bfP}\cdot\left(\begin{matrix}  1 \\0\end{matrix}\right) & = \left(\begin{matrix}{\bfr}_1\\2\lambda_s^2|{\bfr}_1|^2\end{matrix}\right),
\end{split}
\end{equation*}
which proves the claim since $\lambda_s\neq \lambda_2,\lambda_3$. As a result,
 there exists $\varepsilon_1>0$ such that if $|\bfN(\bfomega_s,\bfxi_s,c_s)|<\varepsilon_1$, the entire matrix  $\nabla F$ is regular at $(\lambda_s, \omega_s,\xi_s)$. However, from what we noticed regarding the dependence of $\bfM_s$ on $a$ and  (\ref{boundonN}), this property can be achieved   by taking $a$ sufficiently large, say,  $a>a_1$. Now, since the solutions to (\ref{readytoimplicit}) are characterized by the eigenvalues of the tensor $\bfI(\rho_s)$, by what we have just proved it follows that if the eigenvalues are all distinct and $a$ is large enough, the system (\ref{adieu}) admits only three solution each of which may belong to the $\Omega-$limit set. However, these solutions are isolated and $\Omega$ is connected, and we conclude that $\Omega$ must coincide with one and only one solution. 
We shall next prove that a sufficient condition for  the three eigenvalues of $\bfI(\rho_s)$ to be simple (and, therefore distinct) is that the three eigenvalues of  $\bfI(\bar{\rho})$ be simple  as well, provided $a$ is sufficiently large. To show this, we observe that,
by elementary properties, it is enough to show that there is a sufficiently small $\varepsilon_2$ such that $| \bfI(\bar\rho) - \bfI(\rho_s)|<\varepsilon_2$. We now see that, by definition and (\ref{adieu})$_1$,  $\bfI(\bar\rho) = \bfI(\rho_s)|_{\bfomega_s= {\bf 0}, c_s= (m_\F/|\Ce|)^{\gamma - 1}}$. As a result, we have
$$
| \bfI(\bar\rho) - \bfI(\rho_s)|  \leq \left|\int_0^1 \frac{\partial \bfI}{\partial \sigma}(\sigma \bfomega_s, c_s) \, {\rm d}\sigma\right| + \left|\int_{(m_\F/|\Ce|)^{\gamma -1}}^{c_s} \frac{\partial \bfI}{\partial c}\,({\bf 0}, c) {\rm d}c\right| \equiv \mathcal I_1+\mathcal I_2
$$
Arguing as in the proof of (\ref{boundonN}), we show
\begin{equation}\label{i1}
\mathcal I_1\leq c_1 \,a^{-1/(\gamma-1)}\,.
\end{equation}
As far as  $\mathcal I_2$, since $\frac{\partial \bfI}{\partial c}({\bf 0},c)$ is bounded, we have
$$
\mathcal I_2 \leq c_2 \left|c_s - (m_\F/|\Ce|)^{\gamma - 1}\right|.
$$ 
Equations (\ref{adieu})$_{1,5}$ in combination with the implicit function theorem yields a smooth $c_s=c_s(\eta,\bfomega_s,\bfx_s,m_F)$ with $\eta:=a^{-1/(\gamma-1)}$ and such that $c_s(0,\bfomega_s,\bfxi_s,m\F)=(m_\F/|\Ce|)^{\gamma-1}$. We thus conclude, for $\eta$ small enough
\begin{equation*}
\left|c_s - (m_\F/|\Ce|)^{\gamma - 1}\right|\leq c_3\eta
\end{equation*}
which concludes the proof of the lemma.
\QED\end{proof}

Before stating our main result, we would like to make some comments about the physical meaning of the tensor $\bfI(\bar\rho)$. To this end, suppose we replace the compressible fluid, $\mathcal F$, in the cavity with a fluid, $\bar{\mathcal F}$, of constant density $\bar{\rho}\equiv m_{\mathcal F}/|\mathcal C|$. Also, denote by $\bar{G}$ the center of mass of the coupled system $\bar{\mathcal S}:=\mathcal B\cup\bar{\mathcal F}$. Then (Lemma \ref{stab}) $\bfI(\bar{\rho})$ is the inertia tensor of $\bar{\mathcal S}$ with respect to $\bar G$.
\smallskip\par
The findings of this section, which constitute the main achievement of this paper,  are now collected in the following.
\Bt\label{final.claim}
Let $\Ce$ be of class $C^4$ and let $(\bfu,\rho)$ be a generic solution given in Theorem \ref{long.time.strong}. Suppose that the three eigenvalues of the tensor $\bfI(\bar{\rho})$  are all distinct. Then, there exists $a_0>0$ such that if $a>a_0$,   $(\rho,\bfu)$ tends, as $t\to\infty$, in appropriate norms to a uniquely determined solution $(\rho_s,\bfomega_s,\bfxi_s)$ to (\ref{adieu}). Therefore, the terminal motion of the coupled system $\mathcal S$ reduces to a uniform rotation around an axis parallel to the (constant) angular momentum, $\bfM_0$, of $\mathcal S$ and passing through its center of mass $G$. 
\Et
\begin{proof} In view of what we have already proved, we should only  comment about the statement regarding the terminal motion. However, it is clear that, since in such a state $\bfv\equiv {\bf 0}$ and $\bfomega$ is time-independent, the system $\mathcal S$ rotates as a unique rigid body with constant angular velocity. Moreover, the axis of rotation, ${\sf a}$, must be parallel to $\bfM_0$  by the conservation of angular momentum. Finally, ${\sf a}$ passes through $G$, because $G$ is at rest in any possible motion of $\mathcal S$ and so, in particular, in the terminal one. 

\QED\end{proof}
{\bf Acknowledgments}. Part of this work was carried out when G.P. Galdi was tenured with  the Eduard
\v{C}{e}ch Distinguished Professorship at the Mathematical Institute of the Czech Academy of Sciences in Prague. His work is also partially
supported by NSF Grant DMS-1614011, and the Mathematical Institute of the Czech
Academy of Sciences (RVO 67985840). The research of V. M\'acha  is supported by GA\v{C}R project P201-16-032308
 and RVO 67985840, and that of \v{S}. Ne\v{c}asov\'a by GA\v{C}R project P201-16-032308 and 
RVO 67985840.

\address{galdi@pitt.edu\\ Department of Mechanical Engineering\\ and Materials Science\\
University of Pittsburgh\smallskip\\ and\smallskip\\ macha@math.cas.cz,\ matus@math.cas.cz\\ Institute of Mathematics\\ Czech Academy of Sciences}

\end{document}